\documentclass[a4paper,10pt,twoside]{article}

\usepackage[utf8]{inputenc}
\usepackage{amsmath, amssymb, amsthm}
\usepackage{amstext, amsfonts, a4}
\usepackage{dsfont}
\usepackage{latexsym}
\usepackage{mathrsfs}
\usepackage[all,cmtip]{xy}
\usepackage{graphicx}
\usepackage{enumerate}
\usepackage{hyperref,color}
\usepackage{stmaryrd}
\usepackage[shortlabels]{enumitem}
\usepackage{tikz-cd}
\tikzcdset{crossing over}

\usepackage{url}
\usepackage[pdftex,a4paper,left=3cm,right=3cm,top=3cm,bottom=3cm]{geometry}

\def\ra{\rightarrow}

\def\lra{\longrightarrow}

\def\lmapsto{\longmapsto}

  \def\sL{\mathscr{L}}
  
 \def\sS{\mathscr{S}} 
\def\sV{\mathscr{V}}

\def\bbA{\mathbb{A}}\def\bbC{\mathbb{C}}
\def\bbF{\mathbb{F}}\def\bbG{\mathbb{G}}

\def\bbN{\mathbb{N}}\def\bbP{\mathbb{P}}
\def\bbQ{\mathbb{Q}}\def\bbT{\mathbb{T}}

\def\bbZ{\mathbb{Z}}

\def\cA{\mathcal{A}}\def\cC{\mathcal{C}}\def\cD{\mathcal{D}}
\def\cH{\mathcal{H}}
\def\cI{I(\overline{F}/ F)}

\def\cM{\mathcal{M}}
\def\cR{\mathcal{R}}

\def\whT{\widehat{T}}

\DeclareMathOperator{\diag}{diag}

\DeclareMathOperator{\End}{End}

\DeclareMathOperator{\Hom}{Hom}
\DeclareMathOperator{\hw}{hw}

\DeclareMathOperator{\Id}{Id}

\DeclareMathOperator{\ind}{ind}

\DeclareMathOperator{\Mod}{Mod}

\DeclareMathOperator{\pr}{pr}

\DeclareMathOperator{\Proj}{Proj}

\DeclareMathOperator{\QCoh}{QCoh}

\DeclareMathOperator{\SingDiag}{SingDiag}

\DeclareMathOperator{\Spec}{Spec}
\DeclareMathOperator{\Max}{Max}

\DeclareMathOperator{\ASph}{Sph}

\DeclareMathOperator{\soc}{soc}

\DeclareMathOperator{\Sym}{Sym}

\DeclareMathOperator{\unr}{unr}

\newtheorem{counter}[subsection]{$\!\!$}
\newenvironment{Def}{\begin{counter} {\bf Definition.}}{\end{counter}}

\newenvironment{Prop}{\begin{counter} {\bf Proposition.}}{\end{counter}}
\newenvironment{Lem}{\begin{counter} {\bf Lemma.}}{\end{counter}}
\newenvironment{Cor}{\begin{counter} {\bf Corollary.}}{\end{counter}}
\newenvironment{Th}{\begin{counter} {\bf Theorem.}}{\end{counter}}

\newenvironment{Pt}{\begin{counter} \rm}{\end{counter}}

\newtheorem{counter*}[subsubsection]{$\!\!$}
\newenvironment{Def*}{\begin{counter*} {\bf Definition.}}{\end{counter*}}
\newenvironment{Not*}{\begin{counter*} \rm {\bf Notation.}}{\end{counter*}}
\newenvironment{Notss*}{\begin{counter*} \rm {\bf Notations.}}{\end{counter*}}
\newenvironment{DefNot*}{\begin{counter*} \rm {\bf Definition-Notation.}}{\end{counter*}}
\newenvironment{Nots*}{\begin{counter*} \rm {\bf Notations.}}{\end{counter*}}
\newenvironment{Prop*}{\begin{counter*} {\bf Proposition.}}{\end{counter*}}
\newenvironment{Lem*}{\begin{counter*} {\bf Lemma.}}{\end{counter*}}
\newenvironment{Cor*}{\begin{counter*} {\bf Corollary.}}{\end{counter*}}
\newenvironment{Th*}{\begin{counter*} {\bf Theorem.}}{\end{counter*}}
\newenvironment{Rem*}{\begin{counter*} \rm {\bf Remark.}}{\end{counter*}}
\newenvironment{Ex*}{\begin{counter*} \rm {\bf Example.}}{\end{counter*}}
\newenvironment{Exs*}{\begin{counter*} \rm {\bf Examples.}}{\end{counter*}}
\newenvironment{Pt*}{\begin{counter*} \rm}{\end{counter*}}
\newenvironment{Q*}{\begin{counter*} \rm {\bf Question.}}{\end{counter*}}

\newcommand{\iso}{\stackrel{\sim}{\longrightarrow}}

\title{\textbf{\huge{Lubin-Tate moduli space of semisimple mod $p$ Galois representations for $GL_2$ \\ and Hecke modules}}}

\author{Cédric Pépin and Tobias Schmidt}
\date{\today}

\begin{document}

\maketitle

\begin{abstract}
Let $p$ be an odd prime. Let $F$ be a non-archimedean local field of residue characteristic $p$, and let $\bbF_q$ be its residue field. Let $\cH^{(1)}_{\bbF_q}$ be the pro-$p$-Iwahori-Hecke algebra of the $p$-adic group ${\rm GL_2}(F)$ with coefficients in $\bbF_q$, and let $Z(\cH^{(1)}_{\bbF_q})$ be its center. We define a scheme $X(q)_{\bbF_q}$ whose geometric points parametrize the semisimple two-dimensional Galois representations of ${\rm Gal}(\overline{F}/F)$ over $\overline{\bbF}_q$. Then we construct a morphism from $\Spec Z(\cH^{(1)}_{\bbF_q})$ to $X(q)_{\bbF_q}$ generalizing the morphism appearing in \cite{PS2} for $F=\bbQ_p$. In the case $F/\bbQ_p$, we show that the induced map from Hecke modules to Galois representations, when restricted to supersingular modules, coincides with Grosse-Kl\"onne's bijection \cite{GK18}. For this, we determine the Lubin-Tate $(\varphi,\Gamma)$-modules associated to absolutely irreducible Galois representations. 
\end{abstract}

\tableofcontents

\section{Introduction}

Let $p$ be an odd prime and $F$ a non-archimedean complete local field with ring of integers $o_F$ and residue field $\bbF_q$ of characteristic $p$. Let $\overline{F}$ be an algebraic closure of $F$ and denote by $\overline{\bbF}_q$ its residue field. Let
${\rm Gal}(\overline{F}/F)$ be the absolute Galois group of $F$. Let $\cH^{(1)}_{\bbF_q}$ be the pro-$p$ Iwahori-Hecke algebra of the $p$-adic group
${\rm GL_2}(F)$ with coefficients in $\bbF_q$, and let $Z(\cH^{(1)}_{\bbF_q})$ be its center. When $F=\bbQ_p$, we constructed in \cite{PS2} a 
 morphism 
$$\sL: \Spec Z(\cH^{(1)}_{\bbF_p})\longrightarrow X$$
to the moduli scheme $X$ of semisimple two-dimensional mod $p$ representations of ${\rm Gal}(\overline{\bbQ}_p/\bbQ_p)$ introduced by Emerton-Gee in \cite{Em19}, with the following property: the push-forward along $\sL$ of the extended mod $p$ spherical module $\cM_{\overline{\bbF}_p}^{(1)}$ realizes the semisimple mod $p$ Langlands correspondence for ${\rm GL_2}(\bbQ_p)$ defined by Breuil \cite{Be11}.
\vskip5pt
In the present work, we construct a Lubin-Tate version of the morphism $\sL$ for the local field $F$, with the property that it induces the correspondence defined by Grosse-Kl\"onne when $F/\bbQ_p$, cf. \cite{GK18}.
\vskip5pt
We start by defining a certain two-dimensional $\bbF_p$-scheme $X(q)$ depending only on the parameter $q$. Its connected components are families of chains of projective lines, parametrized by the multiplicative group $\bbG_m$, and $X(q)$ coincides with $X$ above when $q=p$. Our first main result (Thm. \ref{moduli}) is that the geometric points of $X(q)$ parametrize the isomorphism classes of semisimple two-dimensional mod $p$ representations of ${\rm Gal}(\overline{F}/F)$ over $\overline{\bbF}_q$. Writing $q=p^f$, the parametrization depends on the Lubin-Tate fundamental character 
$\omega_f: {\rm Gal}(\overline{F}/F)\rightarrow\bbF_q^\times$ associated with a choice of uniformizer $\pi\in o_F$. 

\vskip5pt 
To go further, we will denote by ${\bf GL_2}$ the Langlands dual group of ${\rm GL}_2$ over the coefficient field $\bbF_q$ and by $W$ its Weyl group. Recall the extended Vinberg toric variety
$V^{(1)}_{\mathbf{\whT}}\rightarrow\bbA^{1}$ associated with the diagonal torus $\mathbf{\whT}\subset {\bf GL_2}$ and its special fibre at $0\in\bbA^1$
$$V^{(1)}_{\mathbf{\whT},0}=\mathbf{\whT}(\bbF_q)\times \SingDiag_{2\times2}\times\bbG_m.$$
%$$
The geometry of the center $Z(\cH^{(1)}_{\bbF_q})$ is best understood in terms of the 
mod $p$ pro-$p$ Iwahori-Satake isomorphism 
\cite[Thm.B]{PS} $$
\xymatrix{
\sS^{(1)}_{\bbF_q}: \Spec Z(\cH^{(1)}_{\bbF_q})\ar[r]^<<<<<{\sim} & S(q).
}
$$
Here, the Satake scheme $S(q):=V^{(1)}_{\mathbf{\whT},0}/W$ is the quotient of $V^{(1)}_{\mathbf{\whT},0}$modulo its natural $W$-action. We show  that there is a completely natural quotient morphism of $\bbF_q$-schemes
$$L: S(q)\longrightarrow X(q)_{\bbF_q}$$ 
to the base change $X(q)_{\bbF_q}$ of $X(q)$ (Thm. \ref{Langlands_morphism}). Its construction is elementary algebraic geometry and does not make use of the Galois parametrization of $X(q)$. 
For example, on generic (regular) connected components of $S(q)$, the morphism $L$ is just the toric construction of the projective line (times $\bbG_m$). In a second step, we precompose the morphism $L$ with the isomorphism $\sS^{(1)}_{\bbF_q}$ 
to obtain a morphism $$\sL:\Spec Z(\cH^{(1)}_{\bbF_q})\longrightarrow X(q)_{\bbF_q}.$$
It gives back the morphism $\sL$ appearing in \cite{PS2} when $F=\bbQ_p$. 
In general, the morphism $\sL$ satisfies several compatibilities, e.g. with regard to twist by characters or Serre weights, which we discuss in sections $8$ and $9$. 
\vskip5pt 
Next, recall the extended mod $p$ spherical module $\cM^{(1)}_{\overline{\bbF}_q}$ from \cite[7.4.1]{PS}. It is a distinguished $\cH^{(1)}_{\overline{\bbF}_q}$-action on the maximal 
commutative subring $\cA^{(1)}_{\overline{\bbF}_q}$ of $\cH^{(1)}_{\overline{\bbF}_q}$ and a mod $p$ analogue (plus extension to the pro-$p$ Iwahori level) of the classical spherical module appearing in complex Kazhdan-Lusztig theory \cite[3.9]{KL87}. The quasi-coherent module (associated to)
$\cM^{(1)}_{\overline{\bbF}_q}$ over $\Spec Z(\cH^{(1)}_{\overline{\bbF}_q})$, when specialized at the closed points of $\Spec Z(\cH^{(1)}_{\overline{\bbF}_q})$, can be used to obtain a parametrization of all irreducible $\cH^{(1)}_{\overline{\bbF}_q}$-modules \cite[7.4.9/7.4.15]{PS}. Combined with the morphism $\sL$, we get a correspondence parametrized by closed points
\[
\begin{tikzcd}[row sep=scriptsize, column sep=scriptsize]
 & &z\in \Max Z(\cH^{(1)}_{\overline{\bbF}_q})\arrow[dr] \arrow[dl] \\
&(\cM^{(1)}_{\overline{\bbF}_q})_z& & \rho_{\sL(z)} \\
\end{tikzcd}
\]
between the $\cH^{(1)}_{\overline{\bbF}_q}$-modules $(\cM^{(1)}_{\overline{\bbF}_q})_z$ and 
the semisimple Galois representations $\rho_{\sL(z)}: {\rm Gal}(\overline{F}/F)\rightarrow {\rm GL_2}(\overline{\bbF}_q)$. When $F/\bbQ_p$, we show  that the singular locus of this correspondence is 1-1 between supersingular irreducible $\cH^{(1)}_{\overline{\bbF}_q}$-modules and irreducible Galois representations; more precisely, we show that it agrees with the bijection established by Grosse-Kl\"onne \cite{GK18} in the case of $GL_2(F)$ (Thm. \ref{comparison_GK}).

\vskip5pt
The construction in \cite{GK18} goes through mod $p$ Lubin-Tate 
$(\varphi,\Gamma)$-modules and their relation \cite{KR09,Sch17} to 
mod $p$ representations of ${\rm Gal}(\overline{F}/F)$. So for our comparison with \cite{GK18},
it is necessary to classify the Lubin-Tate $(\varphi,\Gamma)$-modules corresponding to the absolutely irreducible mod $p$ representations of ${\rm Gal}(\overline{F}/F)$. 
In the cyclotomic case $F=\bbQ_p$, this is a result of Berger \cite{Be10}. We adapt Berger's proof to the Lubin-Tate setting. As in his case, there is no point in restricting to two-dimensional modules and we obtain our classification in any dimension  (Thm. \ref{thm_main}).
\vskip5pt
We recall some background on Lubin-Tate $(\varphi,\Gamma)$-modules in an appendix. 

\vskip5pt

Notation: We keep the notation from the introduction. Let $p>2$ be an odd prime. $F$ denotes a non-archimedean complete local field with ring of integers $o_F$ and residue field of characteristic $p$ and cardinal $q=p^f$. We fix an algebraic closure $\overline{F}/F$, denote by 
$\overline{\bbF}_q/\bbF_q$ its residue field extension, and by $\bbF_{q^n}\subset \overline{\bbF}_q$ the unique subextension of cardinality $q^n$, for each $n\geq 1$. 

For $n\geq 1$, we will denote by ${\rm GL}_n$ the reductive group scheme of invertible $n\times n$-matrices over $F$, and use the same notation for its canonical model over $o_F$ and its special fiber over $\bbF_q$. We will denote by ${\bf GL_n}$ the Langlands dual group of ${\rm GL}_n$ over the coefficient field $\bbF_q$, and use the same notation for its base change to $\overline{\bbF}_q$. 

All Galois representations are supposed to be continuous. 
\vskip5pt
The second author thanks Laurent Berger for answering some questions on $(\varphi,\Gamma)$-modules.

\section{Some reminders on mod $p$ Galois representations}\label{sec_repgalois}

We recall some facts and fix some notation on mod $p$ Galois representations. 
\begin{Pt} 
Let $\pi\in o_F$ be a uniformizer. For any integer $n\geq 1$, let 
$\pi_{nf}\in \overline{F}$ be an element such that $\pi_{nf}^{q^n-1}=-\pi$. We then have 
Serre's fundamental character of level $nf$
$$\omega_{nf}: {\rm Gal}(\overline{F}/F_{n}) \longrightarrow \bbF_{q^n}^{\times}$$ given by 
$ g\mapsto \frac{g(\pi_{nf})}{\pi_{nf}} \in \mu_{q^n-1}(\overline{F})$ followed by reduction mod $\pi$, cf. \cite{Se72}.
One has $$\omega_{nf}^{\frac{q^n-1}{q-1}}=\omega_f|_{{\rm Gal}(\overline{F}/F_{n})}.$$

Let $\cI\subset {\rm Gal}(\overline{F}/F)$ be the inertia subgroup and let 
$\cI^t$ be its tame quotient. Choose an element $\varphi\in {\rm Gal}(\overline{F}/F)$ lifting the Frobenius $x\mapsto x^q$ on $ {\rm Gal}(\overline{\bbF}_q/\bbF_q)$. Since $\omega_{f}: \cI \rightarrow\bbF^{\times}_q$ is surjective \cite[Prop. 2]{Se72}, we may and will assume $\omega_{f}(\varphi)=1$. Note that the restriction $$\omega_{nf}: \cI \longrightarrow \bbF_{q^n}^{\times}$$
of the character $\omega_{nf}$ to $\cI$ is canonical, since we defined $\overline{\bbF}_q$ as the residue field of $\overline{F}$ and 
$\bbF_{q^n}^{\times}$ as its unique subfield of cardinality $q^n$.

\end{Pt}
\begin{Pt}\label{local_class_field}
We normalize local class field theory $F^{\times}\rightarrow {\rm Gal}(\overline{F}/F)^{\rm ab}$ by sending 
$\pi$ to the geometric Frobenius $\varphi^{-1}$. 
In this way, we identify the smooth $\overline{\bbF}^{\times}_q$-valued characters of ${\rm Gal}(\overline{F}/F)$ and of $F^{\times}$.  This restricts to a bijection between smooth characters of the inertia subgroup $\cI$ and of $\bbF_q^\times$.
For example, the fundamental character $\omega_f: F^{\times} \ra \overline{\bbF}^{\times}_q$
is the extension of the inclusion $\omega: \bbF_q^{\times}\stackrel{\subset} {\longrightarrow}\overline{\bbF}^{\times}_q$ to $F^{\times}$ 
satisfying $\omega_f(\pi)=1$. 
\end{Pt}

\begin{Pt} Let $F\subseteq F_n\subseteq \overline{F}$ be the unique unramified extension of degree $n$ over $F$. A smooth character $\chi: {\rm Gal}(\overline{F}/F_n)\rightarrow\overline{\bbF}^\times_q$ is {\it regular} if its ${\rm Gal}(F_n/F)$- conjugates $\chi, \chi^{q},..., \chi^{q^{n-1}}$ are all distinct. The irreducible smooth $\overline{\bbF}_q$-representations of ${\rm Gal}(\overline{F}/F) $ of dimension $n$ are 
given by the representations 
$$ \ind^{ {\rm Gal}(\overline{F}/F) }_{ {\rm Gal}(\overline{F}/F_n) }(\chi)$$
smoothly induced from the regular $\overline{\bbF}_q$-characters $\chi$ of ${\rm Gal}(\overline{F}/F_n).$ The conjugates $\chi, \chi^{q},..., \chi^{q^{n-1}}$ of $\chi$ induce isomorphic representations and there are no other isomorphisms between the representations \cite[1.14]{V94}, \cite[5.1]{V04}.
\end{Pt}

\begin{Pt}
A character $\omega_{nf}^h$ for $1\leq h \leq q^n-2$ is regular if and only if its conjugates $\omega_{nf}^{h}, \omega_{nf}^{qh},..., \omega_{nf}^{q^{n-1}h}$ are all distinct. Equivalently, if and only if $h$ is $q$-primitive, that is, there is no $d<n$
such that $h$ is a multiple of $(q^n-1)/(q^d-1)$.
 The representation $\ind^{{\rm Gal}(\overline{F}/F)}_{ {\rm Gal}(\overline{F}/F_n) }(\omega_{nf}^h)$ is then defined over $\bbF_{q^n}$. It has a basis 
 $\{ v_0,...,v_{n-1}\}$ of eigenvectors for the characters 
 $\omega_{nf}^{h}, \omega_{nf}^{qh},..., \omega_{nf}^{q^{n-1}h}$ 
 of ${\rm Gal}(\overline{F}/F_n)$ such that $\varphi(v_i)=v_{i-1}$ and $\varphi(v_0)=v_{n-1}$. In particular, its determinant coincides with 
 $\omega^h_f$ on the subgroup
 ${\rm Gal}(\overline{F}/F_{n})$ and takes $\varphi$ to $(-1)^{n-1}$.
  \end{Pt}
  \begin{Pt}
For $\lambda\in\overline{\bbF}_q^{\times}$, let $\mu_{\lambda}$ or $\unr(\lambda)$ be the unramified character of ${\rm Gal}(\overline{F}/F)$ sending $\varphi^{-1}$ to 
$\lambda$. Fix $\delta$ with $\delta^n=(-1)^{n-1}$.  The representation
 $$ \ind (\omega_{nf}^h):= ( \ind^{{\rm Gal}(\overline{F}/F) }_{ {\rm Gal}(\overline{F}/F_n) }(\omega_{nf}^h) ) \otimes\mu_{\delta}$$
 is then uniquely determined by the two conditions 
 $$\det  \ind (\omega_{nf}^h) =  \omega^h_f\hskip10pt \text{and} \hskip10pt\ind (\omega_{nf}^h) |_{\cI} = \omega_{nf}^{h} \oplus \omega_{nf}^{qh} \oplus ...\oplus \omega_{nf}^{q^{n-1}h} .$$ 
 \end{Pt}
\begin{Pt}\label{irred_k_reps}
 Let $\bbF_q\subset k \subset \overline{\bbF}_q$ be an intermediate extension of $\bbF_q$. Every absolutely irreducible smooth $k$-representation of ${\rm Gal}(\overline{F}/F) $ of dimension $n$ is isomorphic to 
$\ind (\omega_{nf}^h) \otimes \mu_{\lambda}$ for a $q$-primitive 
$1\leq h \leq q^n-2$ and a scalar $\lambda\in\overline{\bbF}_q^{\times}$ such that $\lambda^n\in k^{\times}$ and one has 
$$\ind (\omega_{nf}^h) \otimes \mu_{\lambda}\simeq \ind (\omega_{nf}^{\tilde{h}}) \otimes \mu_{\tilde{\lambda}}$$
if and only if  ${\rm Gal}(F_n/F).\omega_{nf}^h={\rm Gal}(F_n/F).\omega_{nf}^{\tilde{h}}$
and $\lambda^n=\tilde{\lambda}^n$.

\vskip5pt 

Since $\omega_{nf}^{\frac{q^n-1}{q-1}}=\omega_f$, every irreducible representation of ${\rm Gal}(\overline{F}/F) $ of dimension $n$ is therefore isomorphic to 
$\ind (\omega_{nf}^h) \otimes \omega_f^s \mu_{\lambda}$ for
a $q$-primitive $1 \leq h\leq \frac{q^n-1}{q-1} -1$, a scalar $\lambda\in \overline{\bbF}^{\times}_q$, and 
$0\leq s \leq q-2$.
 \end{Pt}
 \begin{Pt} \label{irred_n2} Let $n=2$. Since $\frac{q^2-1}{q-1} -1=q$ and $\ind (\omega_{2f})\simeq \ind (\omega_{2f}^q)$, 
every irreducible representation of ${\rm Gal}(\overline{F}/F) $ of dimension $2$ is isomorphic to 
$\ind (\omega_{2f}^h) \otimes \omega_f^s \mu_{\lambda}$ for 
a $q$-primitive $1 \leq h\leq q-1$,  a scalar $\lambda\in \overline{\bbF}^{\times}_q$ and $0\leq s \leq q-2$.  
 \end{Pt}

\section{The $q$-parametrization of two-dimensional inertial types}

A \emph{tame (two-dimensional) mod $p$ inertial type} is (the isomorphism class of) a 
continuous homomorphism $$\tau : I(\overline{F}/ F)^t\rightarrow {\bf GL_2}(\overline{\bbF}_q),$$
which extends to a representation of ${\rm Gal}(\overline{F}/F)$. 
A tame mod $p$ inertial type is semi-simple as a representation of $I(\overline{F}/ F)^t$, so we can write $\tau=\chi_1\oplus\chi_2$. In the terminology of \cite[sec.11]{H09}, $\tau$ is of \emph{niveau $1$} if $\chi_i^{q-1}=1$ for $i=1,2$ and of \emph{niveau $2$} otherwise. 

\subsection{The basic even case}\label{types_basic_even_case}

We will give a parametrization of the tame mod $p$ inertial types with determinant $\omega_f$.
\begin{Pt*}
For $r\in\{0,2,\ldots,q-5,q-3\}$, set $s(r):=-\frac{r}{2}\in\{0,-1,\ldots,-\frac{q-5}{2},-\frac{q-3}{2}\}$. Then consider the table with two columns:
$$
\begin{array}{ccc}
(0,0)  &  | & (q-3,1)  \\
(2,-1) & |  &  (q-5,2) \\
 \vdots &  \vdots & \vdots   \\
(r,s(r)) & |  &  (q-3-r,s(r)+r+1) \\
 \vdots &  \vdots & \vdots   \\
(q-5,-\frac{q-5}{2}) & | & (2,\frac{q-3}{2}) \\
(q-3,-\frac{q-3}{2}) & | &  (0,\frac{q-1}{2}).
\end{array}
$$
\end{Pt*}

\begin{Pt*}
To each pair $(r,s(r))$, 
we attach the type of niveau $1$ with determinant $\omega_f$
$$ \tau:=\left(\begin{array}{cc}
\omega_{f}^{r+1} & 0\\
0 & 1
\end{array} \right)\otimes\omega_{f}^{s(r)}\simeq \left(\begin{array}{cc}
\omega_{f}^{q-2-r} & 0\\
0 & 1
\end{array} \right)
\otimes\omega_{f}^{s(r)+r+1}.
$$
According to the above table, this gives $\frac{q-1}{2}$ types of niveau $1$. 

To each pair $(r,s(r))$, we may also attach a type of niveau $2$ 
with determinant $\omega_f$, namely
$$ \tau:=\left(\begin{array}{cc}
\omega_{2f}^{r+1} & 0\\
0 & \omega_{2f}^{q(r+1)}
\end{array} \right)\otimes\omega_{f}^{s(r)}\simeq \left(\begin{array}{cc}
\omega_{2f}^{q-r} & 0\\
0 & \omega_{2f}^{q(q-r)} 
\end{array} \right)
\otimes\omega_{2f}^{s(r-2)+r-1}.
$$
According to the above table, this gives $\frac{q-1}{2}+1=\frac{q+1}{2}$ types of niveau $2$. 
As in the case of $F=\bbQ_p$ \cite[3.3]{PS2}, one shows that all types 
with determinant $\omega_f$ are obtained in this way. 
\end{Pt*}

\subsection{The basic odd case}\label{types_basic_odd_case}

We will give a parametrization of the tame mod $p$ 
inertial types with determinant $1$.
\begin{Pt*}\label{table_even}

For $r\in\{-1,1,\ldots,q-4,q-2\}$, set $s(r):=-\frac{r+1}{2}\in\{0,-1,\ldots,-\frac{q-3}{2},-\frac{q-1}{2}\}$. Then consider the table with two columns: 
$$
\begin{array}{ccc}
  (-1,0) &  | & (q-2,0) \\
 (1,-1) & |  &(q-4,1)  \\
 \vdots &  \vdots & \vdots   \\
 (r,s(r)) & |  &  (q-3-r,s(r)+r+1)\\
 \vdots &  \vdots & \vdots   \\
 (q-4,-\frac{q-3}{2})& | & (1,\frac{q-3}{2}) \\
 (q-2,-\frac{q-1}{2})& | &(-1,\frac{q-1}{2})  .
\end{array}
$$
\end{Pt*}

\begin{Pt*}
To each pair $(r,s(r))$, 
we attach the type of niveau $1$ with determinant $1$
$$ \tau:=\left(\begin{array}{cc}
\omega_{f}^{r+1} & 0\\
0 & 1
\end{array} \right)\otimes\omega_{f}^{s(r)}\simeq \left(\begin{array}{cc}
\omega_{f}^{q-2-r} & 0\\
0 & 1
\end{array} \right)
\otimes\omega_{f}^{s(r)+r+1}.
$$
According to the above table, this gives $\frac{q+1}{2}$ types of niveau $1$. 

To each pair $(r,s(r))$, we may also attach a type of niveau $2$ 
with determinant $1$, namely
$$ \tau:=\left(\begin{array}{cc}
\omega_{2f}^{r+1} & 0\\
0 & \omega_{2f}^{q(r+1)}
\end{array} \right)\otimes\omega_{f}^{s(r)}\simeq \left(\begin{array}{cc}
\omega_{2f}^{q-r} & 0\\
0 & \omega_{2f}^{q(q-r)} 
\end{array} \right)
\otimes\omega_{2f}^{s(r-2)+r-1}.
$$
According to the above table, this gives $\frac{q-1}{2}$ types of niveau $2$. 
As in the case of $F=\bbQ_p$ \cite[3.3]{PS2}, one shows that all types 
with determinant $1$ are obtained in this way. 
\end{Pt*}

\subsection{The general case}\label{types_general_case}

\begin{Pt*} Let  $$N_q:=\{0,1,\cdots,q-2\} \hskip20pt
 \text{and}\hskip20pt D_q:=1+N_q=\{1,\cdots,q-1\}.$$
It splits into
$$
N_q=E_q\coprod O_q
$$
where
$$
E_q:=\{2m\}_{m=0,\ldots, \frac{q-3}{2}} \hskip20pt \text{resp.} \hskip20pt  
O_q=\{2m+1\}_{m=0,\ldots, \frac{q-3}{2}}
$$
is the subset of even resp. odd numbers, which both have cardinality $(q-1)/2$.
For each $n\in N_q$ let $d_n:=1+n\in D_q$.
\end{Pt*}
\begin{Lem*} \label{number_types_niv1}
Let $n\in N_q$. The number of tame types of niveau $1$ with determinant 
$\omega_{f}^{d_n}$ is
$$
\left \{\begin{array}{cc}
\frac{q-1}{2} & \textrm{if $n\in E_q$} \\
\frac{q+1}{2} & \textrm{if $n\in O_q$}.
\end{array} \right.
$$
In particular, the total number of tame types of niveau $1$ is
$\frac{q^2-q}{2}.$
\end{Lem*}

\begin{proof}
The first part follows from the niveau $1$ part of the 
basic even case (if $n$ even) or basic odd case (if $n$ odd)
by twisting with powers of $\omega_f$. 
The second part follows from this, since 
$$ |E_q|\frac{q-1}{2}+|O_q|\frac{q+1}{2}=\frac{q-1}{2}\Big(\frac{q-1}{2}+\frac{q+1}{2}\Big)=\frac{q^2-q}{2}. $$
\end{proof}
\begin{Lem*} \label{number_types_niv2}
Let $n\in N_q$. The number of tame types of niveau $2$ with determinant 
$\omega_{f}^{d_n}$ is
$$
\left \{\begin{array}{cc}
\frac{q+1}{2} & \textrm{if $n\in E_q$} \\
\frac{q-1}{2} & \textrm{if $n\in O_q$}.
\end{array} \right.
$$
In particular, the total number of tame types of niveau $2$ is
$\frac{q^2-q}{2}.$
\end{Lem*}

\begin{proof}
The first part follows from the niveau $2$ part of the 
basic even case (if $n$ even) or basic odd case (if $n$ odd)
by twisting with powers of $\omega_f$. 
The second part follows from this, since 
$$ |E_q|\frac{q+1}{2}+|O_q|\frac{q-1}{2}=
\frac{q-1}{2}\Big(\frac{q+1}{2}+\frac{q-1}{2}\Big)=\frac{q^2-q}{2}. $$
\end{proof}

\section{The $q$-scheme of semisimple Galois representations}

The following is inspired by the work of Emerton-Gee in the case $F=\bbQ_p$ \cite{Em19}.

\begin{Pt} \textbf{The projective line.}
Let 
$$
\bbP^1:=\Proj(\bbF_p[x,y])
$$
be the projective line over $\bbF_p$. It is the gluing of the two affine lines
$$
\bbA^1_x:=\Spec(\bbF_p[x])\subset \bbP^1 \supset \bbA^1_y:=\Spec(\bbF_p[y])
$$
along the open $\Spec(\bbF_p[x^{\pm1}])=\Spec(\bbF_p[y^{\mp1}])$. The closed complement of $\bbA^1_x$ is the  \emph{the point at infinity} 
$\infty :=[1:0]\in\bbA_y^1(\bbF_p)$ and the closed complement of $\bbA^1_y$ is \emph{the origin} $0 :=[0:1]\in\bbA^1_x(\bbF_p)$:
$$
\bbP^1=\bbA^1_x\cup\{\infty\}=\{0\}\cup\bbA^1_y.
$$
There is a natural $\bbG_m$-action on $\bbP^1$ by "scaling", given by $(\alpha,z):=\alpha z$ for $\alpha\in\bbG_m$ and $z\in\bbA^1_x$. The space of orbits $\bbP^1/\bbG_m$ has three elements,
the two closed orbits $\{ 0,\infty\}$ and the open orbit $\bbP^1\setminus \{ 0,\infty\}$. The action depends on the choice of the affine coordinate $x$, but $\bbP^1/\bbG_m$ does not.
  \end{Pt}

\begin{Pt} \textbf{The chains of $\bbP^1$'s.}
Let $l\in\bbN_{\geq 1}$. For each $i\in\{0,\ldots,l-1\}$, set $\cC_i:=\Proj(\bbF_p[x_i,y_i])$, a copy of the projective line. Then let
$$
\widetilde{\cC}(l):=\cC_0 \coprod \cC_1\coprod  \cdots\coprod \cC_{l-2}\coprod\cC_{l-1}
$$
be the disjoint union of these $l$ copies of $\bbP^1$. Finally let
$$
\cC(l):=\cC_0\; {}_{\infty}\bigcup{}_0\;\cC_1 \;{}_{\infty}\bigcup{}_0 \cdots {}_{\infty}\bigcup{}_0\;\cC_{l-2}\; {}_{\infty}\bigcup{}_0\;\cC_{l-1}
$$
be the $\bbF_p$-scheme obtained by identifying the point at infinity of $\cC_i$ with the origin of $\cC_{i+1}$ for all $i$ from $0$ to $l-2$, a chain of length $l$ of copies of the projective line. In particular, it is a curve over $\bbF_p$, with $l$ irreducible components and $l-1$ singularities which are ordinary simple nodes, and the canonical morphism
$$
\widetilde{\cC}(l)\lra \cC(l)
$$
is its normalization. \end{Pt}

\begin{Pt} \textbf{The even and the odd $q$-chains.}\label{q-chains}
We call
$$
\cC(\frac{q-1}{2})
$$
the \emph{even $q$-chain}. 
There is a natural $\bbG_m$-action on $\cC(\frac{q-1}{2})$ induced by the scaling action on each component $\cC_i$ and we denote by $\cC(\frac{q-1}{2})/\bbG_m=\cup_{i=0,...,\frac{q-3}{2}} \;\cC_i/\bbG_m$ its space of $\bbG_m$-orbits. 

\vskip5pt

We call
$$
\cC(\frac{q+1}{2})
$$
the \emph{odd $q$-chain}. There is a natural $\bbG_m$-action on each "interior" component 
$\cC_i$ for $0<i<\frac{q-1}{2}$. On the smooth part of the two "exterior" components $\bbA^1_x\subset \cC_0$ and 
 $\bbA^1_y\subset \cC_{\frac{q-1}{2}}$, we pretend\footnote{Let $\bbG_m$ act on itself by multiplication. 
 There is no $\bbG_m$-action on $\bbA^1$ making the map $t\mapsto t+t^{-1}$ equivariant.} to have a "modified action" via the parametrization $\bbG_m\rightarrow\bbA^1, t\mapsto t+t^{-1}$. In other words, we {\it define} the 
 space $\cC_0/\bbG_m$ to consist of two elements, the point $\infty$ and its open complement. Similarly, we {\it define} the space $\cC_{\frac{q-1}{2}}/\bbG_m$ to consist of two elements, the origin $0$ and its open complement. 
 Finally, we let 
 $\cC(\frac{q-1}{2})/\bbG_m=\cup_{i=0,...,\frac{q-1}{2}}\;\cC_i/\bbG_m$.
 \end{Pt}

\begin{Pt} \label{cc_X} \textbf{Connected components.}
For each $n\in N_q$ and $d_n=1+n\in D_q$, we define the $\bbF_p$-scheme
$$
X_{d_n}(q):= 
\left \{\begin{array}{cc}
\cC(\frac{q-1}{2}) \times \bbG_m& \textrm{if $n\in E_q$} \\
\cC(\frac{q+1}{2}) \times \bbG_m & \textrm{if $n\in O_q$}.
\end{array} \right.
$$
There is a natural diagonal $\bbG_m$-action on $X_{d_n}(q)$ induced by the $\bbG_m$-action
on $\cC(\frac{q\mp 1}{2})$ (using the "modified action" from \ref{q-chains} on the exterior components of $\cC(\frac{q+1}{2})$) and by
the multiplication with the {\it square } $(\alpha,z)\mapsto \alpha^2z$ on $\bbG_m$\footnote{In the subsequent parametrization of $X_{d_n}(q)$ by semisimple Galois representations \ref{EGparam}, an element $z\in\bbG_m$ will correspond to the \emph{determinant} of the $\varphi$-action, which gets multiplied by the \emph{square} of $\alpha$ when the representation is twisted by $\unr(\alpha)$.}. We let $X_{d_n}(q)/ \bbG_m$ denote the space of $\bbG_m$-orbits.
\end{Pt}

\begin{Def*}\label{moduli_scheme} We define the \emph{$q$-scheme of semisimple two-dimensional mod $p$ Galois representations} to be the scheme over $\bbF_p$ 
$$
X(q):=\coprod_{n\in N_q} X_{d_n}(q).
$$
There is a natural diagonal $\bbG_m$-action on $X(q)$ coming from the $\bbG_m$-action
on each connected component $X_{d_n}(q)$. We let $X(q)/ \bbG_m$ denote the space of $\bbG_m$-orbits.
\end{Def*}
The terminology for $X(q)$ will become clear in the next subsection. 
\vskip5pt
\textbf{Remark.} Fix once and for all a generator
$\zeta$ of $\bbF_q^{\times}$.  The $\bbG_m$-action on $X(q)$ can in fact be promoted to an action of the group $\bbG_m(\bbF_q)\times \bbG_m$: the action of the element $(\zeta^n,z)\in\bbG_m(\bbF_q)\times\bbG_m$ on the connected component $X_{d_m}(q)$ is given by the map 

$$
\xymatrix{
\Id\times\{z^2\}:X_{d_m}(q) \ar[r]^{\;\;\;\;\;\sim} & X_{d_{m+2n}}(q)}.
$$
Here, we take the representative of $m+2n$ modulo $(q-1)$ in $N_q$ and $\{z^2\}$ refers to the automorphism of $\bbG_m$ given by multiplication by $z^2$.

\begin{Pt*}
Let $D(q)$ be the finite constant $\bbF_p$-scheme such that $D(q)(\bbF_p)=D_q$.
The scheme $X(q)$ is canonically fibered over $D(q)$: 
$$
d(q):X(q)\lra D(q)\quad \textrm{ with $d(q)^{-1}(d_n)=X_{d_n}(q)$ for all $d_n\in D_q$}.
$$
It also admits a canonical projection $\pr_2$ to $\bbG_m$, whence a canonical morphism
$$
d(q)\times\pr_2:X(q)\lra D(q)\times\bbG_m.
$$
From now on, we drop the $(q)$ from the notation, so we write $X$ instead of $X(q)$ and so on.
\end{Pt*}

\begin{Pt}\textbf{The $q$-Galois parametrization.}\label{EGparam}
Recall that we have fixed an arithmetic Frobenius $\varphi\in{\rm Gal}(\overline{F}/ F)$. 
In the preceding subsection we have defined a certain $\bbF_p$-scheme $X$.
The aim of the present subsection is to establish the following theorem. 
\end{Pt}

\begin{Th*}\label{moduli}
There is a canonical (up to a sign) bijection $$
\iota_{\varphi}:X(\overline{\bbF}_q) \cong \big\{\text{semisimple $\rho : {\rm Gal}(\overline{F}/ F)\rightarrow \mathbf{\bf GL_2}(\overline{\bbF}_q)$} \big\}/ \sim.
$$

More precisely, for each $n\in N_q$ and $z_2\in \bbG_m(\overline{\bbF}_q)$, there is a canonical (up to a sign) bijection:
$$
\iota_{\varphi,n}|_{\pr_2=z_2}:X_{d_n}|_{\pr_2=z_2}(\overline{\bbF}_q) 
\cong \big\{\text{semisimple $\rho : {\rm Gal}(\overline{F}/ F)\rightarrow {\bf GL_2}(\overline{\bbF}_q)\ |\ \det (\rho)= \omega_{f}^{d_n}\unr(z_2)$} \big\}/ \sim.
$$
\end{Th*}

The proof proceeds along the lines of the case $F=\bbQ_p$ \cite{Em19}, cf. also \cite{PS2}, 
by assigning to the geometric standard coordinates $(x,y),z_2$ on each irreducible component $\bbP^1\times\bbG_m$ of $X$ an isomorphism class of semisimple representations 
$\rho:{\rm Gal}(\overline{F}/F)\ra {\bf GL_2}(\overline{\bbF}_q)$. 

\begin{Pt*}\label{basic_even_case} \textbf{The basic even case.} Let us consider the case where
$$
n=0\in E_q\quad\textrm{i.e. $d_n=1\in D_q$},\quad\textrm{and}\quad z_2=1\in\bbG_m.
$$
Then
$$
X_{d_0}|_{\pr_2=1}=
\cC_0\; {}_{\infty}\bigcup{}_0\;\cC_1 \;{}_{\infty}\bigcup{}_0 \cdots {}_{\infty}\bigcup{}_0\;\cC_{\frac{q-5}{2}}\; {}_{\infty}\bigcup{}_0\;\cC_{\frac{q-3}{2}}\times\{1\}.
$$
For $i\in\{0,1,\ldots,\frac{q-5}{2},\frac{q-3}{2}\}$, set $r:=2i$. Then $r\in\{0,2,\ldots,q-5,q-3\}$, and we rewrite the above chain as
$$
X_{d_0}|_{\pr_2=1}=
\cC^0\; {}_{\infty}\bigcup{}_0\;\cC^2\; {}_{\infty}\bigcup{}_0 \cdots {}_{\infty}\bigcup{}_0\;\cC^{q-5} \;{}_{\infty}\bigcup{}_0\;\cC^{q-3}\times\{1\}.
$$
Next, for $r\in\{0,2,\ldots,q-5,q-3\}$, set $s(r):=-\frac{r}{2}\in\{0,-1,\ldots,-\frac{q-5}{2},-\frac{q-3}{2}\}$. Then reconsider the table from \ref{types_basic_even_case}: 

$$
\begin{array}{ccc}
(0,0)  &  | & (q-3,1)  \\
(2,-1) & |  &  (q-5,2) \\
 \vdots &  \vdots & \vdots   \\
(r,s(r)) & |  &  (q-3-r,s(r)+r+1) \\
 \vdots &  \vdots & \vdots   \\
(q-5,-\frac{q-5}{2}) & | & (2,\frac{q-3}{2}) \\
(q-3,-\frac{q-3}{2}) & | &  (0,\frac{q-1}{2}).
\end{array}
$$
It gives a rule to attach an isomorphism class of representations $\rho$ to a point of $[x:1]=[1:y]\in\bbP^1\setminus\{0,\infty\}=\cC^r\setminus\{0,\infty\}$: one takes
$$
\rho:=\left(\begin{array}{cc}
\unr(x)\omega_{f}^{r+1} & 0\\
0 & \unr(x^{-1})
\end{array} \right)
\otimes\omega_{f}^{s(r)}
\simeq
\left(\begin{array}{cc}
\unr(y)\omega_{f}^{q-2-r} & 0\\
0 & \unr(y^{-1})
\end{array} \right)
\otimes\omega_{f}^{s(r)+r+1}.
$$
Moreover, one takes
$$
\rho:=\ind(\omega_{2f}^{r+1})\otimes\omega_{f}^{s(r)}
$$
at the origin $0$, and
$$
\rho:=\ind(\omega_{2f}^{q-3-r})\otimes\omega_{f}^{s(r)+r+1}\simeq
\ind(\omega_{2f}^{r+3})\otimes\omega_{f}^{s(r+2)}
$$
at the point $\infty$. We have thus a well-defined map 
$$
\iota_{\varphi,0}|_{\pr_2=1}: X_{d_0}|_{\pr_2=1}(\overline{\bbF}_q) \rightarrow \big\{\text{semisimple continuous $\rho : {\rm Gal}(\overline{F}/ F)\rightarrow {\bf GL_2}(\overline{\bbF}_q)\ |\ \det (\rho )=\omega_{f}$} \big\}/ \sim.
$$

By its very construction, it is compatible with the parametrization of mod $p$ tame inertial types in the basic even case, cf. \ref{types_basic_even_case}.
\end{Pt*}

\begin{Lem*}\label{Comp} The map $\iota_{\varphi,0}|_{\pr_2=1}$ induces bijections
$$
\text{open orbits in  $X_{d_0}/\bbG_m$}
\simeq \big\{\text{types $\tau$ of niveau $1$ $|$ $\det \tau= \omega_{f}$} \big\}
$$
$$
\text{closed orbits in  $X_{d_0}/\bbG_m$}
\simeq \big\{\text{types $\tau$ of niveau $2$ $|$ $\det \tau= \omega_{f}$} \big\}.
$$
\end{Lem*}

\begin{proof} This follows from the case $n=0$ in \ref{number_types_niv1} and
\ref{number_types_niv2}.
\end{proof}

\begin{Cor*} 
The map $\iota_{\varphi,0}|_{\pr_2=1}$ is a bijection. 
\end{Cor*}

\begin{Pt*} \label{basic_odd_case}\textbf{The basic odd case.} Let us consider the case where
$$
n=q-2\in O_q\quad\textrm{i.e. $d_n=q-1\in D_q$},\quad\textrm{and}\quad z_2=1\in\bbG_m.
$$
Then
$$
X_{d_{q-2}}|_{\pr_2=1}=
\cC_0\; {}_{\infty}\bigcup{}_0\;\cC_1 {}_{\infty}\bigcup{}_0 \cdots {}_{\infty}\bigcup{}_0\;\cC_{\frac{q-3}{2}}\; {}_{\infty}\bigcup{}_0\;\cC_{\frac{q-1}{2}}\times\{1\}.
$$
For $i\in\{0,1,\ldots,\frac{q-3}{2},\frac{q-1}{2}\}$, set $r:=2i-1$. Then $r\in\{-1,1,\ldots,q-4,q-2\}$, and we rewrite the above chain as
$$
X_{d_{q-2}}|_{\pr_2=1}=
\cC^{-1}\; {}_{\infty}\bigcup{}_0\;\cC^1\; {}_{\infty}\bigcup{}_0 \cdots {}_{\infty}\bigcup{}_0\;\cC^{q-4}\; {}_{\infty}\bigcup{}_0\;\cC^{q-2}\times\{1\}.
$$
Next, for $r\in\{-1,1,\ldots,q-4,q-2\}$, set $s(r):=-\frac{r+1}{2}\in\{0,-1,\ldots,-\frac{q-3}{2},-\frac{q-1}{2}\}$. Then reconsider the table from \ref{types_basic_odd_case}:

$$
\begin{array}{ccc}
  (-1,0) &  | & (q-2,0) \\
 (1,-1) & |  &(q-4,1)  \\
 \vdots &  \vdots & \vdots   \\
 (r,s(r)) & |  &  (q-3-r,s(r)+r+1)\\
 \vdots &  \vdots & \vdots   \\
 (q-4,-\frac{q-3}{2})& | & (1,\frac{q-3}{2}) \\
 (q-2,-\frac{q-1}{2})& | &(-1,\frac{q-1}{2})  .
\end{array}
$$
For $r\notin\{-1,q-2\}$, it gives a rule to attach an isomorphism class of representations $\rho$ to a point of $[x:1]=[1:y]\in\bbP^1\setminus\{0,\infty\}=\cC^r\setminus\{0,\infty\}$: one takes
$$
\rho:=\left(\begin{array}{cc}
\unr(x)\omega_{f}^{r+1} & 0\\
0 & \unr(x^{-1})
\end{array} \right)
\otimes\omega_{f}^{s(r)}
\simeq
\left(\begin{array}{cc}
\unr(y)\omega_{f}^{q-2-r} & 0\\
0 & \unr(y^{-1})
\end{array} \right)
\otimes\omega_{f}^{s(r)+r+1}.
$$
Moreover, one takes
$$
\rho:=\ind(\omega_{2f}^{r+1})\otimes\omega_{f}^{s(r)}
$$
at the origin $0$, and
$$
\rho:=\ind(\omega_{2f}^{q-3-r})\otimes\omega_{f}^{s(r)+r+1}\simeq
\ind(\omega_{2f}^{r+3})\otimes\omega_{f}^{s(r+2)}
$$
at the point $\infty$. For $r\in\{-1,q-2\}$, one uses the surjection
\begin{eqnarray*}
\bbG_m & \lra& \bbA^1 \\
 z & \lmapsto & t:=z+z^{-1}
\end{eqnarray*}
and attach to $t\in\bbA^1=\cC^{-1}\setminus\{\infty\}$, resp. $t\in\bbA^1=\cC^{q-2}\setminus\{0\}$):
$$
\rho:=\left(\begin{array}{cc}
\unr(z) & 0\\
0 & \unr(z^{-1})
\end{array} \right)
\simeq
\left(\begin{array}{cc}
\unr(z^{-1}) & 0\\
0 & \unr(z)
\end{array} \right)
$$
resp.
$$
\rho:=\left(\begin{array}{cc}
\unr(z) & 0\\
0 & \unr(z^{-1})
\end{array} \right)\otimes\omega_{f}^{\frac{q-1}{2}}
\simeq
\left(\begin{array}{cc}
\unr(z^{-1}) & 0\\
0 & \unr(z)
\end{array} \right)\otimes\omega_{f}^{\frac{q-1}{2}}.
$$
We have thus a well-defined map $$\iota_{\varphi,q-2}|_{\pr_2=1}: X_{d_{q-2}}|_{\pr_2=1}(\overline{\bbF}_q) \rightarrow \big\{\text{semisimple continuous $\rho : {\rm Gal}(\overline{F}/ F)\rightarrow {\bf GL_2}(\overline{\bbF}_q)\ |\ \det \rho= 1$} \big\}/ \sim.
$$

By its very construction, it is compatible with the parametrization of mod $p$ tame inertial types in the basic odd case, cf. \ref{types_basic_odd_case}.
\end{Pt*}

\begin{Lem*}\label{Comp} The map $\iota_{\varphi,q-2}|_{\pr_2=1}$ induces bijections
$$
\text{open orbits in  $X_{d_{q-2}}/\bbG_m$}
\simeq \big\{\text{types $\tau$ of niveau $1$ $|$ $\det \tau= 1$} \big\}
$$
$$
\text{closed orbits in  $X_{d_{q-2}}/\bbG_m$}
\simeq \big\{\text{types $\tau$ of niveau $2$ $|$ $\det \tau= 1$} \big\}.
$$
\end{Lem*}

\begin{proof} This follows from the case $n=q-2$ in \ref{number_types_niv1} and
\ref{number_types_niv2}.
\end{proof}

\begin{Cor*} 
The map $\iota_{\varphi,q-2}|_{\pr_2=1}$ is a bijection. 
\end{Cor*}

\begin{Pt*} \label{qparam_the_general_case}\textbf{The general case.} Let
$$
n\in N_q= E_q\coprod O_q\quad\textrm{i.e. $d_n\in D_q$},\quad\textrm{and}\quad z_2\in\bbG_m(\overline{\bbF}_q).
$$
Choose a square root $\sqrt{z_2}$ of $z_2$. This choice is responsable to the addendum {\it up to a sign} in the statement of the theorem \ref{moduli}. Set
$$
\eta:=
\unr(\sqrt{z_2})\left\{\begin{array}{cc}
\omega_{f}^{\frac{d_n-1}{2}} & \textrm{if $n\in E_q$} \\
\omega_{f}^{\frac{d_n}{2}} & \textrm{if $n\in O_q$}.
\end{array}\right.
$$
Then there is obviously a unique (bijective) map 
$$
X_{d_n}|_{\pr_2=z_2}(\overline{\bbF}_q)\xrightarrow{\sim} \big\{\text{semisimple $\rho : {\rm Gal}(\overline{F}/ F)\rightarrow {\bf GL_2}(\overline{\bbF}_q)\ |\ \det \rho = \omega_{f}^{d_n}\unr(z_2)$} \big\}/ \sim
$$
such that the bijections
$$
\xymatrix{
\Id\times\{z_2\}:X_{d_0}|_{\pr_2=1}(\overline{\bbF}_q) \ar[r]^{\sim} & X_{d_n}|_{\pr_2=z_2}(\overline{\bbF}_q)\quad  \textrm{if $n\in E_q$}
}
$$
and
$$
\xymatrix{
\Id\times\{z_2\}:X_{d_{q-2}}|_{\pr_2=1}(\overline{\bbF}_q) \ar[r]^{\sim} & X_{d_n}|_{\pr_2=z_2}(\overline{\bbF}_q)\quad  \textrm{if $n\in O_q$}
}
$$
correspond on the Galois side to twisting by the character $\eta$.
Here $\{z_2\}$ refers to the automorphism of $\bbG_m$ given by multiplication by $z_2$.
This ends the proof of the theorem \ref{moduli}.
\end{Pt*}
As a corollary of the proof, we obtain 

\begin{Cor*}\label{Comp} The map $\iota_{\varphi}$ induces
a bijection 
$$
\text{$X/\bbG_m$}
\simeq \big\{\text{tame inertial types $\tau: \cI\rightarrow {\bf GL_2}(\overline{\bbF}_q)$}\big\}.
$$
Under this bijection open and closed orbits correspond to types of niveau $1$ and $2$ respectively.
\end{Cor*}
\begin{Pt*} \label{qparam_twisting}\textbf{Twisting.} The Galois parametrization in theorem \ref{moduli} is compatible with twisting in the following sense. Recall the action of $\bbG_m(\bbF_q)\times \bbG_m$ on $X$, depending on our fixed generator $\zeta$ of 
$\bbF_q^\times$, cf. remark after definition \ref{moduli_scheme}.
Now the group $\bbG_m(\bbF_q)\times \bbG_m$ is naturally isomorphic to the group of Galois characters via $(\zeta^n,z_2)\mapsto \omega_f^n\unr(z_2)$. Let 
$$\cR_{m,z_2}:=\big\{\text{semisimple $\rho : {\rm Gal}(\overline{F}/ F)\rightarrow {\bf GL_2}(\overline{\bbF}_q)\ |\ \det (\rho)= \omega_{f}^{d_m}\unr(z_2)$} \big\}/ \sim$$
be the set appearing on the right hand-side of theorem \ref{moduli}. Suppose $m\in E_q$ and let $\eta:=\unr(\sqrt{z_2})\omega_{f}^{\frac{d_m-1}{2}}$ be a "choice of sign" inducing the Galois parametrization $$\iota_{\varphi,m}|_{\pr_2=z_2}: X_{d_m}|_{\pr_2=z_2}\simeq  \cR_{m,z_2}.$$
Let an arbitrary Galois character
$\chi:=\omega_f^{n}\unr(z)$ be given.  

It leads to the sign choice 
$\eta':=\unr(\sqrt{z_2}\cdot z)\omega_{f}^{\frac{d_{m+2n}-1}{2}}$ and the Galois parametrization 
$$\iota_{\varphi,m+2n}|_{\pr_2=z_2z^2}: X_{d_{m+2n}}|_{\pr_2=z_2z^2}\simeq  \cR_{m+2n,z_2z^2}.$$
The two Galois parametrizations make the left-hand side and the back-side of the 
following half cube
\[
\begin{tikzcd}[row sep=scriptsize, column sep=scriptsize]
& X_{d_0}|_{\pr_2=1}\arrow[dr] \arrow[rr] \arrow[dd] & & X_{d_m}|_{\pr_2=z_2} \arrow[dl, "{(\zeta^n,z)}", blue] \arrow[dd,blue] \\  & &X_{d_{m+2n}}|_{\pr_2=z_2z^2} \arrow[dd,blue] \\
& \cR_{0,1} \arrow[dr, "\eta' "] \arrow[rr,"\hskip40pt\eta"] & & \cR_{m,z_2}  \arrow[dl, "\chi",blue] \\
& & \cR_{m+2n,z_2z^2}  \\
\end{tikzcd}
\]
commutative. The top of the half cube is commutative, since trivially
$$ \Id\times\{z_2z^2\}=\Id\times\{z^2\}\circ \Id\times\{z_2\}.$$
The bottom of the half cube is commutative, since $\eta'=\eta\cdot\chi$. It follows that the right-hand side (written in blue) is commutative as well. This means that the action of 
$(\zeta^n,z)\in \bbG_m(\bbF_q)\times \bbG_m$ on $X_{d_m}|_{\pr_2=z_2}$ corresponds
to twisting by the corresponding Galois character $\chi=\omega_f^{n}\unr(z)$ on the corresponding set of Galois representations
$\cR_{m,z_2}$. The case $m\in O_q$ is similar.

\end{Pt*}

\section{The $q$-scheme of Satake parameters} 

We recall some notions and results from \cite{PS} and \cite{PS2}. In the following, all schemes and fiber products are over $\Spec \bbF_q$.
Let $W$ be the Weyl group of ${\bf GL_2}$ and $w$ its nontrivial element. 

\begin{Pt} \label{intro_vinbergtorus} 
Let $\mathbf{\whT}$ be the torus of invertible diagonal $2\times2$ matrices. We consider the scheme 
$$
V_{\mathbf{\whT},0} := \SingDiag_{2\times2}\times\;\bbG_m,
$$
where $\SingDiag_{2\times2}$ represents the semigroup of singular diagonal $2\times 2$-matrices \cite[7.1]{PS}. Consider the extended semigroup
$$V^{(1)}_{\mathbf{\whT},0}:=\mathbf{\whT}(\bbF_q)\times V_{\mathbf{\whT},0}.$$ 
It carries a natural $W$-action: the natural action of $W$ on 
the factors $\mathbf{\whT}(\bbF_q)$ and  $\SingDiag_ {2\times 2}$ and the trivial one on $\bbG_m$. 
\end{Pt}

\begin{Def} We define the \emph{scheme of mod $p$ Satake parameters for ${\bf GL_2}$} to be the scheme over $\bbF_q$
$$S(q):= V^{(1)}_{\mathbf{\whT},0}/W.$$
\end{Def}

\begin{Pt}
The scheme $S(q)$ is {\it canonically} fibered over the finite constant $\bbF_q$-scheme $\mathbf{\whT}(\bbF_q)/W$:
$$
\pi_0:S(q)\lra \mathbf{\whT}(\bbF_q)/W.
%\quad \textrm{ with $\pi_0^{-1}(\gamma)=V_{\mathbf{\whT},0}^{\gamma}/W$ for all $\gamma\in \mathbf{\whT}(\bbF_q)/W$}.
$$
The fibers of $\pi_0$ are the connected components of $S(q)$. The irreducible components of 
$S(q)$ can be labelled by the elements of $\mathbf{\whT}(\bbF_q)$. This depends on a choice of order $(t_1,t_2)$ on every regular orbit $\{t_1\neq t_2\}$ in $\mathbf{\whT}(\bbF_q)/W$: the order induces an isomorphism 
$\pi_0^{-1}(\{t_1, t_2\})\simeq V_{\mathbf{\whT},0}=\bbA_x^1\cup_0\bbA_y^1$ and we can label the image of $\bbA^1_x$ (resp. of $\bbA^1_y$) in 
$\pi_0^{-1}(\{t_1, t_2\})$ by $t_1$ (resp. $t_2$).

Composing with the determinant map $\mathbf{\whT}(\bbF_q)/W\ra\bbG_m(\bbF_q)$ gives a morphism 
$$S(q)\lra \mathbf{\whT}(\bbF_q)/W\lra \bbG_m(\bbF_q).$$
\end{Pt}

\begin{Pt}
The scheme $S(q)$ also admits the canonical projection $\pr_2:S(q)\ra\bbG_m$. Whence finally a composed morphism
$$
S(q)\lra \mathbf{\whT}(\bbF_q)/W\times\bbG_m\lra\bbG_m(\bbF_q)\times\bbG_m.
$$

From now on, we drop the $(q)$ from the notation, i.e. we will write $S$ instead of $S(q)$ and so on.
\end{Pt}

\section{The $\bbF_q$-morphism $L$ from Satake to Galois} \label{Fq_Langlands_mor} 
%We keep the previous notation. 
The aim of the present section is to establish the following theorem. Let $X_{\bbF_q}$ be the base change to $\bbF_q$ of the $\bbF_p$-scheme $X=X(q)$ of semisimple two-dimensional Galois representations, cf. \ref{moduli_scheme}.

\begin{Pt}\textbf{The twisting action.} \label{L_twisting}
Recall our choice of generator $\zeta$ of the group $\bbF_q^{\times}$. 
According to the remark after definition \ref{moduli_scheme}, there is a natural action of the group 
$\bbG_m(\bbF_q)\times\bbG_m$ on the scheme 
$X$. This action extends linearly to $X_{\bbF_q}$. On the other hand, define the following action of $\bbG_m(\bbF_q)\times\bbG_m$ on 
$$
S=V^{(1)}_{\mathbf{\whT},0}/W=\big(\mathbf{\whT}(\bbF_q)\times \SingDiag_{2\times2}\times\bbG_m\big)/W.
$$
Let $(n,z_2)\in N_q\times \bbG_m$. The element $\zeta^n\in\bbG_m(\bbF_q)$ acts only via the factor $\mathbf{\whT}(\bbF_q)$, by multiplication by 
$\diag(\zeta^n,\zeta^n)$. The element $z_2\in\bbG_m$ acts trivially on $\mathbf{\whT}(\bbF_q)$, by multiplication by $\diag(z_2,z_2)$ on 
$\SingDiag_{2\times2}$, and by multiplication by $z_2^2$ on $\bbG_m$.
\end{Pt}
\begin{Th}\label{Langlands_morphism}
There is a quotient morphism of $\bbF_q$-schemes 
$$L: S\lra X_{\bbF_q}$$ 
which gives back the morphism $L$ appearing in \cite{PS2} in the case 
$F=\bbQ_p$. The morphism $L$ is $\bbG_m(\bbF_q)\times\bbG_m$-equivariant.
\end{Th} 
As in \cite{PS2}, the morphism is a quotient morphism, locally given by the toric construction of the projective line (except on the two exterior components in the odd case, see below). Its construction goes along the lines of \cite{PS2}. 

For each $(n,z_2)\in N_q\times\bbG_m$, we have the fibre $S_{(\zeta^n,z_2)}$ of $S$ at $(\zeta^n,z_2)\in \bbG_m(\bbF_q)\times\bbG_m$, and the fiber $X_{(d_n,z_2)}$ of $X_{\bbF_q}$  at $(d_n,z_2)\in D_q\times\bbG_m$. We will have a morphism 
$$L_{(n,z_2)}: S_{(\zeta^n,z_2)}\longrightarrow X_{(d_n,z_2)}.$$

\begin{Pt}\textbf{The ordering on the irreducible components of $S$.} \label{ordering_comp_S}
Let $(n,z_2)\in N_q\times\bbG_m$. Let $x,y$ resp. $z_1$ be the canonical standard coordinates resp. Steinberg coordinate on each regular resp. non-regular connected component of $S_{(\zeta^n,z_2)}$. According to \cite[8.1.2]{PS}, we have the following description of $S_{(\zeta^n,z_2)}$.

\vskip5pt

Suppose $n$ is even, i.e. $n\in E_q$. Then $S_{(\zeta^n,z_2)}$ is the disjoint union 
$$
 \bbA^1_{z_1} \;\coprod \;\bbA^1_x\cup_0 \bbA^1_y \;\coprod \cdot\cdot\cdot \coprod \;\bbA^1_x\cup_0 \bbA^1_y\; \coprod  \;\bbA^1_{z_1}\times\{z_2\}
 $$
indexed by the fibre of $\mathbf{\whT}(\bbF_q)/W\ra\bbG_m(\bbF_q)$ in $\zeta^n$. 
%There is a given ordering on this fibre and hence on the elements of the disjoint union: 
The irreducible components of $S_{(\zeta^n,z_2)}$ can be labelled by the sequence of \emph{ordered} pairs of elements of $\mathbf{\whT}(\bbF_q)$
$$
t_i \cdot \diag(\zeta^s,\zeta^s), \hskip5pt t_i^w \cdot \diag(\zeta^s,\zeta^s)
$$ 
where $n=2s$ and $t_i:=\diag(\zeta^i,\zeta^{-i})$ (and $t_i^w$ its $w$-conjugate) for $i=0,...,\frac{q-1}{2}$.  Choose a square root $\sqrt{z_2}$. The action of the element $(\zeta^s,\sqrt{z_2})\in \bbG_m(\bbF_q)\times\bbG_m$ gives an isomorphism
$$(\zeta^s,\sqrt{z_2}): S_{(\zeta^0,1)}\iso S_{(\zeta^n,z_2)}$$
which preserves the ordering.
\vskip5pt
Suppose $n$ is odd, i.e. $n\in O_q$. Then $S_{(\zeta^n,z_2)}$ is the disjoint union 

$$ \bbA^1_x\cup_0 \bbA^1_y \;\coprod \cdot\cdot\cdot \coprod \;\bbA^1_x\cup_0 \bbA^1_y \times\{z_2\}$$
indexed by the fibre of $\mathbf{\whT}(\bbF_q)/W\ra\bbG_m(\bbF_q)$ in $\zeta^n$. 
%There is a given ordering on this fibre and hence on the elements of the disjoint union: 
The irreducible components of $S_{(\zeta^n,z_2)}$ can be labelled by the sequence of \emph{ordered} pairs of elements of $\mathbf{\whT}(\bbF_q)$ 
$$
t_i \cdot \diag(\zeta^s,\zeta^s), \hskip5pt t_i^w \cdot \diag(\zeta^s,\zeta^s)
$$ 
where $n=2s-1$ and $t_i:=\diag(\zeta^{i-1+\frac{q-1}{2}},\zeta^{-i+\frac{q-1}{2}})$ (and $t_i^w$ is its $w$-conjugate) for $i=1,...,\frac{q-1}{2}$. 
%Remark:Here the $s$ is different from in the definition of the set $O_q$, since there we have $n=2m+1$ instead of $n=2s-1$. Is it possible to arrange things so that we can avoid this difference between m and s ?
The action of the element $(\zeta^s,\sqrt{z_2})\in \bbG_m(\bbF_q)\times\bbG_m$ gives an isomorphism
$$ (\zeta^s, \sqrt{z_2}): S_{(\zeta^{q-2},1)}\iso S_{(\zeta^n,z_2)}$$
which preserves the ordering.

\end{Pt}

\begin{Pt*}\textbf{The morphism $L$ in the even case.}\label{L_even_case}
Let  $n\in E_q$. We restrict first to the case $n=0$ and $z_2=1$.
% i.e. $(\zeta^n,z_2)=(1,1)$ and $(d_n,z_2)=(1,1)$. 
We have by definition 
$$ X_{(d_0,1)}=
\cC_{0,\bbF_q}\; {}_{\infty}\bigcup{}_0\;\cC_{1,\bbF_q}\; {}_{\infty}\bigcup{}_0 \cdots {}_{\infty}\bigcup{}_0\;\cC_{\frac{q-5}{2},\bbF_q}\; {}_{\infty}\bigcup{}_0\;\cC_{\frac{q-3}{2},\bbF_q}\times\{1\}.
$$
Let $Q_i$ be the origin $0$ on $\cC_{i,\bbF_q}$ for $i\in\{0,1,\ldots,\frac{q-5}{2},\frac{q-3}{2}\}$ and let $Q_{\frac{q-1}{2}}$ be the point $\infty$ on 
$\cC_{\frac{q-3}{2},\bbF_q}$. On the other hand, let $P_i$ be the origin on the $i$-th connected component of $S_{(\zeta^0,1)}$ for $i=0,...,\frac{q-1}{2}.$ 
We map the sequence of points $P_i$ to the sequence of points $Q_i$, i.e. we define
$$L_{(0,1)} (P_i):=Q_{i}.$$ 
Next, suppose $0<i<\frac{q-1}{2}$ and consider the $i$-th connected component  $\bbA^1_x\cup_0 \bbA^1_y$ of $S_{(\zeta^0,1)}$.
Then $L_{(0,1)}(P_i)=Q_i \in \cC_{i-1,\bbF_q}\cap\cC_{i,\bbF_q}$. We define 
$$ L_{(0,1)}(0,y):=[1:y] \in \cC_{i-1,\bbF_q}\hskip5pt \text{and} \hskip5pt L_{(0,1)}(x,0):=[x:1] \in \cC_{i,\bbF_q}.$$
%Note that the defined maps are $\bbG_m$-equivariant in an obvious sense.
Finally, if $i=0$ resp. $i=\frac{q-1}{2}$ we call the Steinberg variable $z_1$ simply $x$ resp. $y$ and put 
$$L_{(0,1)}(x):=[x:1] \in\cC_{0,\bbF_q}\hskip5pt \text{resp.} \hskip5pt  L_{(0,1)}(y):=[1:y] \in\cC_{\frac{q-3}{2},\bbF_q}.$$
%Again, these maps are $\bbG_m$-equivariant in an obvious sense.
We have defined a quotient morphism of $\bbF_q$-schemes  $$L_{(0,1)}: S_{(\zeta^0,1)}\longrightarrow X_{(d_0,1)}$$
which, locally, is the toric construction of the projective line: it identifies the open subset $\bbG_m$ in the "first" irreducible component $\bbA^1$ of a connected component of $S_{(\zeta^0,1)}$ with the open subset $\bbG_m$ in the "second" irreducible component 
$\bbA^1$ of the "next" connected component via the map $z\mapsto z^{-1}$, thus forming a 
$\bbP^1$.

%Moreover, the morphism 
%$\sL_{(0,1)}$ is $\bbG_m$-equivariant (when we ignore the action on the $\bbG_m$-part of the source and the target). 

Let now $n=2s\in E_q$ and $z_2\in\bbG_m$ be general. 
The action of $(\zeta^s,\sqrt{z_2})\in \bbG_m(\bbF_q)\times\bbG_m$ gives
the isomorphism $$(\zeta^s,\sqrt{z_2}): X_{(d_0,1)}\iso X_{(d_n,z_2)}.$$ We define $L_{(d_n,z_2)}:=(\zeta^s,\sqrt{z_2}) \circ L_{(0,1)}  \circ(\zeta^s,\sqrt{z_2}) ^{-1}$. It is well-defined, i.e. independent of the choice of square root $\sqrt{z_2}$.

\end{Pt*}

\begin{Pt*}\textbf{The morphism $L$ in the odd case.} \label{L_odd_case} 
Let  $n\in O_q$. We restrict first to the case $n=q-2$ and $z_2=1$. We have by definition 
$$ X_{(d_{q-2},1)}=
\cC_{0,\bbF_q}\; {}_{\infty}\bigcup{}_0\;\cC_{1,\bbF_q}\; {}_{\infty}\bigcup{}_0 \cdots {}_{\infty}\bigcup{}_0\;\cC_{\frac{q-3}{2},\bbF_q}\; {}_{\infty}\bigcup{}_0\;\cC_{\frac{q-1}{2},\bbF_q}\times\{1\}.
$$
Moreover, on $\cC_{0,\bbF_q}$ we write $t$ for the variable $x$ (so that the double point is at $t=\infty$) and
on $\cC_{\frac{q-1}{2},\bbF_q}$ we write $t$ for the variable $y$ (so that the double point is at $t=\infty$, again). 

Now let $Q_i$ be the origin $0$ on $\cC_{i,\bbF_q}$ for $i\in\{1,2,\ldots,\frac{q-3}{2},\frac{q-1}{2}\}$. On the other hand, let $P_i$ be the origin on the $i$-th connected component of $S_{(\zeta^{q-2},1)}$ for $i=1,...,\frac{q-1}{2}.$ 
We map the sequence of points $P_i$ to the sequence of points $Q_{i}$, i.e. we define
$$L_{(q-2,1)} (P_i):=Q_i.$$ 
Next consider the $i$-th connected component  $\bbA^1_x\cup_0 \bbA^1_y$ of $S_{(\zeta^{q-2},1)}$.
Then $L_{(q-2,1)}(P_i)=Q_i \in \cC_{i-1,\bbF_q}\cap\cC_{i,\bbF_q}$. We define 
$$ L_{(q-2,1)}(0,y):=[1:y] \in \cC_{i-1,\bbF_q}\hskip5pt\text{(if $i\neq 1$)} \hskip15pt \text{and} \hskip15pt L_{(q-2,1)}(x,0):=[x:1] \in \cC_{i,\bbF_q} \hskip5pt\text{(if $i\neq \frac{q-1}{2}$)} .$$
%Note that the defined maps are $\bbG_m$-equivariant in an obvious sense.
Finally, if $i=1$ resp. $i=\frac{q-1}{2}$ we call $t$ the standard variable $y$ resp. $x$ on the $i$-th connected component of $S_{(\zeta^{q-2},1)}$ and put 
$$L_{(q-2,1)}(0,t):=t+t^{-1} \in\cC_{0,\bbF_q}\hskip5pt \text{resp.} \hskip5pt  L_{(q-2,1)}(t,0):=t+t^{-1}\in\cC_{\frac{q-1}{2},\bbF_q}.$$
%By definition of the "modified action" on the target, these maps are "$\bbG_m$-equivariant". 
We have defined a quotient morphism of $\bbF_q$-schemes  $$L_{(q-2,1)}: S_{(\zeta^{q-2},1)}\longrightarrow X_{(d_{q-2},1)}$$
which, locally, is the toric construction of the projective line, except on the two "outer" irreducible components $\bbA^1$ of $S_{(\zeta^{q-2},1)}$, where it is the covering $\bbA^1\rightarrow\bbP^1, t\mapsto t+t^{-1}$.
%Moreover, the morphism 
%$\sL_{(q-2,1)}$ is $\bbG_m$-equivariant (when we ignore the action on the $\bbG_m$-part of the source and the target). 

Let now $n=2s-1\in O_q$ and $z_2\in\bbG_m$ be general. The action of $(\zeta^s,\sqrt{z_2})\in \bbG_m(\bbF_q)\times\bbG_m$ gives the isomorphism $$(\zeta^s,\sqrt{z_2}): X_{(d_{q-2},1)}\iso X_{(d_n,z_2)}.$$ We define $L_{(d_n,z_2)}:=(\zeta^s,\sqrt{z_2}) \circ L_{(q-2,1)}  \circ(\zeta^s,\sqrt{z_2}) ^{-1}$. It is well-defined, i.e. independent of the choice of square root $\sqrt{z_2}$.
\end{Pt*}

\section{The $\bbF_q$-morphism $\sL$ from Hecke to Galois} 

\begin{Pt}\label{Notation_Hecke}
Let $\cH^{(1)}_{\bbF_q}$ be the pro-$p$ Iwahori-Hecke algebra of the group
${\rm GL_2}(F)$ with coefficients in $\bbF_q$, and let $Z(\cH^{(1)}_{\bbF_q})$ be its center. 

Set $\bbT:=T(\bbF_q)$, and denote by $\bbT^{\vee}$ its group of characters. The group $W$ acts naturally on $\bbT$ and $\bbT^{\vee}$, and the connected components of the scheme $\Spec(Z(\cH^{(1)}_{\bbF_q}))$ are \emph{canonically} indexed by the quotient set $\bbT^{\vee}/W$. Moreover, set ${\small
u=
\left (\begin{array}{cc}
0 & 1\\
\pi & 0
\end{array} \right)
}$
and $U=T_u\in \cH^{(1)}_{\bbF_q}$\footnote{This element corresponds to what is denoted by $u^{-1}$ in \cite{PS,PS2}.}. Then $U^2$ is a free invertible element of $Z(\cH^{(1)}_{\bbF_q})$. Whence a canonical morphism of $\bbF_q$-schemes
$$
\pi_0\times\pr_{\Spec(\bbF_q[U^{\pm2}])}:\Spec(Z(\cH^{(1)}_{\bbF_q}))\lra \bbT^{\vee}/W\times\Spec(\bbF_q[U^{\pm2}]).
$$
Finally, restricting along the diagonal cocharacter $\bbF_q^{\times}\ra \bbT$ induces a map $\bbT^{\vee}/W\ra(\bbF_q^{\times})^{\vee}$, whence a composed morphism 
$$
\Spec(Z(\cH^{(1)}_{\bbF_q}))\lra \bbT^{\vee}/W\times\Spec(\bbF_q[U^{\pm2}])\lra (\bbF_q^{\times})^{\vee}\times\Spec(\bbF_q[U^{\pm2}]).
$$
\end{Pt}

\begin{Pt} 
In \cite[Thm.B]{PS} we established the mod $p$ pro-$p$-Iwahori Satake isomorphism
$$
\xymatrix{
\sS^{(1)}_{\bbF_q}: \Spec Z(\cH^{(1)}_{\bbF_q})\ar[r]^<<<<<{\sim} & S(q).
}
$$

Recall our fixed choice of generator $\zeta$ of $\bbF_q^{\times}$. Using the evaluation of cocharacters of $\mathbf{\whT}$ on $\zeta$, one gets an identification of $\bbT^{\vee}$ with $\mathbf{\whT}(\bbF_q)$. Similarly, recall our fixed choice of inclusion $\bbF_q \subset \overline{\bbF}_q$ inducing the character $\omega : \bbF_q^{\times}\ra \overline{\bbF}_q^{\times}$; one gets an identification of $(\bbF_q^{\times})^{\vee}=\langle \omega \rangle$ with 
$\bbG_m(\bbF_q)=\langle \zeta \rangle$.

Then, by construction, the isomorphism $\sS^{(1)}_{\bbF_q}$ fits into the commutative diagram
$$
\xymatrix{
\Spec Z(\cH^{(1)}_{\bbF_q})\ar[r]^<<<<<<<<<<{\sS^{(1)}_{\bbF_q}}_{\sim} \ar[d]_{\pi_0\times\pr_{\Spec(\bbF_q[U^{\pm2}])}} & S(q) \ar[d]^{\pi_0\times \pr_2}\\
\bbT^{\vee}/W\times\Spec(\bbF_q[U^{\pm2}]) \ar@{=}[r] \ar[d] & \mathbf{\whT}(\bbF_q)/W\times\bbG_m \ar[d] \\
(\bbF_q^{\times})^{\vee}\times\Spec(\bbF_q[U^{\pm2}]) \ar@{=}[r] & \bbG_m(\bbF_q)\times\bbG_m.
}
$$
\end{Pt}

\begin{Def}
Composing $\sS^{(1)}_{\bbF_q}$ with the morphism $L: S(q)\ra X_{\bbF_q}$
from \ref{Langlands_morphism} yields a \emph{Langlands morphism}
$$
\sL:= L\circ \sS^{(1)}_{\bbF_q}: \Spec Z(\cH^{(1)}_{\bbF_q})\lra X_{\bbF_q}.
$$
\end{Def}

\begin{Pt}\label{twistHecke} 
Viewing a $\cH^{(1)}_{\bbF_q}$-module as a quasi-coherent module on 
$\Spec Z(\cH^{(1)}_{\bbF_q})$ yields the functor 
$$\sL_*: \Mod (\cH^{(1)}_{\bbF_q} ) \lra \QCoh(X_{\bbF_q}),$$
generalizing the one from \cite[7.2]{PS2} in the case of $F=\bbQ_p$.
% and fixed central character 

Furthermore, let $\cM^{(1)}_{\bbF_q}$ be the mod $p$ spherical module, cf. 
\cite[Thm. A]{PS}. Tensoring $\cM^{(1)}_{\bbF_q}$ over 
$Z(\cH^{(1)}_{\bbF_q})$ with an $\overline{\bbF}_q$-valued central character defines the spherical map 
 $$
\xymatrix{
\ASph:  \big(\Spec Z(\cH^{(1)}_{\bbF_q})\big)(\overline{\bbF}_q)\ar[r] & \{\textrm{left $\cH^{(1)}_{\overline{\bbF}_q}$-modules}\}/\sim.
}
$$
It induces a parametrization of {\it all} irreducible $\cH^{(1)}_{\overline{\bbF}_q}$-modules. In general, modules of the form  $\ASph(v)$ are of length one or two, and they are always of length one if $v^*:Z(\cH^{(1)}_{\bbF_q})\ra \overline{\bbF}_q$ is a supersingular central character, cf. \cite[Thm. E]{PS}. 
\end{Pt}

\begin{Prop} \label{thm_LT} 
The morphism $\sL$ induces a bijection 
 $$
\left.
\begin{array}{lll}
\big(\Spec Z(\cH^{(1)}_{\bbF_q})\big)(\overline{\bbF}_q)^{\rm supersing} & \stackrel{\sim} {\longrightarrow}X_{\bbF_q}(\overline{\bbF}_q)^{\rm irred}& 
\end{array}
\right.
$$
between the sets of \emph{supersingular simple} Hecke modules, via $\ASph$, and of \emph{irreducible} Galois representations, via $\iota_{\varphi}$.
\end{Prop}

\begin{proof}
Let $(n,z_2)\in N_q\times \bbG_m$. The isomorphism $\sS^{(1)}_{\bbF_q}$ induces an isomorphism between the fibers $(\Spec Z(\cH^{(1)}_{\bbF_q}))_{(\omega^n,z_2)}$ and $S_{(\zeta^n,z_2)}$, which in turn is mapped onto $X_{(d_n,z_2)}$ by $L$. The supersingular central characters in 
$(\Spec Z(\cH^{(1)}_{\bbF_q})_{(\omega^n,z_2)}$ correspond to the points $(x,y)=0$ (in the regular case) and $z_1=0$ (in the non-regular case) in $S_{(\zeta^n,z_2)}$. Moreover, by construction of $L_{(n,z_2)}$ and $\iota_{\varphi}$, these points are mapped in a $1:1$ way to points in $X_{(d_n,z_2)}(\overline{\bbF}_q)$ corresponding to irreducible Galois representations. Letting $n$ vary, the resulting injective map 
 $$
\left.
\begin{array}{lll}
\big(\Spec Z(\cH^{(1)}_{\bbF_q})\big)_{z_2}(\overline{\bbF}_q)^{\rm supersing} & \stackrel{\sim} {\longrightarrow}X_{z_2}(\overline{\bbF}_q)^{\rm irred}& 
\end{array}
\right.
$$
is bijective, since source and target have the same cardinality $\frac{q^2-q}{2}$, cf. \cite[Rem. 5.1]{V04}.
\end{proof}

\section{Relation to Grosse-Kl\"onne's functor} 

Combining the spherical map $\ASph$ with 
the morphism $\sL$ gives a correspondence 
$$ \ASph(v)\rightsquigarrow \rho_{\sL(v)}$$
from (certain) $\cH^{(1)}_{\overline{\bbF}_q}$-modules to semisimple Galois representations
${\rm Gal}(\overline{F}/F)\ra {\bf GL_2}(\overline{\bbF}_q)$.
In \cite{PS2} we have shown that in the case $F=\bbQ_p$, the correspondence
$ \ASph(v)\rightsquigarrow \rho_{\sL(v)}$
{\it is} the semisimple mod $p$ local Langlands correspondence\footnote{The category of smooth mod $p$ representations of ${\rm GL_2}(\bbQ_p)$ is equivalent to the category of $\cH^{(1)}_{\overline{\bbF}_q}$-modules \cite{O09}.}
 for the group ${\rm GL_2}(\bbQ_p)$. 

 In the general case $F/\bbQ_p$, Proposition \ref{thm_LT} shows that $$ \ASph(v)\rightsquigarrow \rho_{\sL(v)}$$ induces a bijection between simple supersingular Hecke modules and irreducible Galois representations. In this section, we will show that this bijection 
 {\it is} the (functorial) bijection in the case $n=2$ constructed by Grosse-Kl\"onne \cite{GK18}.
 This makes use of the case $n=2$ in the classification of  irreducible
 \'etale mod $p$ Lubin-Tate $(\varphi,\Gamma)$-modules, the main result \ref{thm_main} of the appendix.
 
 \begin{Pt} \label{Fq_action} 
 To start with, recall the twisting action of $\bbG_m(\bbF_q)\times\bbG_m$ on $S$ from 
 \ref{L_twisting}. Under the isomorphism $\sS^{(1)}_{\bbF_q}$, it corresponds to an action of $ (\bbF_q^{\times})^{\vee}\times\bbG_m$ on
$\Spec(Z(\cH_{\bbF_q}^{(1)}))$, which, in particular, induces an action of $(\bbF_q^{\times})^{\vee}$ on $\bbT^{\vee}/W$. In the sequel, we will denote the later as
 $$  \bbT^{\vee}/W \times (\bbF_q^{\times})^{\vee} \longrightarrow  \bbT^{\vee}/W, \;(\gamma, \omega^n) \mapsto \gamma.\omega^n.$$
\end{Pt}

\begin{Pt}\label{twist}
Let
$$
\xymatrix{
\ASph^{\rm ss}:  \big(\Spec Z(\cH^{(1)}_{\bbF_q})\big)(\overline{\bbF}_q)\ar[r] & \{\textrm{semisimple left $\cH^{(1)}_{\overline{\bbF}_q}$-modules}\}/\sim
}
$$
be the composition of $\ASph$ followed by semisimplification. It is then equivariant for the action of $(\bbF_q^{\times})^{\vee}\times \bbG_m$ on the target deduced from the following twisting action of irreducible (or, more generally, standard) $\cH^{(1)}_{\overline{\bbF}_q}$-modules, cf. \cite[2.5]{PS2}. Let $(n,z_2)\in N_q\times \bbG_m$. In the non-regular case, the $U$-action gets multiplied by $z_2$, the $S$-action remains unchanged and the component $\gamma$ gets multiplied by $\omega^n$ as above. 
%cf. \cite[1.6]{V04}. 
In the regular case, the actions of $X,Y,U^2$ get multiplied by $z_2,z_2,z_2^{2}$ respectively and the component $\gamma$ gets multiplied by 
$\omega^n$ again.
% cf. \cite[2.4]{V04}. 
\end{Pt}

\begin{Pt}\label{ZGveeZvee}
Let $Z(G)$ be the center of $G:={\rm GL}_2(F)$. It is isomorphic to $F^{\times}$ via the diagonal cocharacter $F^{\times}\ra {\rm GL_2}(F)$. Denote by $Z(G)^{\vee}$ the group of smooth $\overline{\bbF}_q$-valued characters of $Z(G)$. It is isomorphic to $(\bbF_q^{\times})^{\vee}\times \bbG_m$ via
$$
Z(G)^\vee\iso(\bbF_q^{\times})^{\vee}\times \bbG_m,\; \;\; \eta\mapsto (\eta |_{\bbF_q^{\times}},\eta(\pi)).
$$
Indeed, any smooth character $F^{\times}\ra \overline{\bbF}_q^{\times}$ is trivial on the subgroup $1+\pi o_F$, and we have normalized local class field theory by sending $\pi$ to $\varphi^{-1}$ in \ref{local_class_field}. 
\end{Pt}

\begin{Pt}
From \ref{twist} and \ref{ZGveeZvee}, we get a twisting action of $Z(G)^\vee$ on standard $\cH^{(1)}_{\overline{\bbF}_q}$-modules, that we denote by
$$
(M,\eta)\lmapsto M\otimes\eta.
$$
Let $\bbF_q\subseteq k \subset \overline{\bbF}_q$ be a finite extension. Let $\cH^{(1)}_k$ be the pro-$p$ Iwahori-Hecke algebra with coefficients in $k$. %In particular, $\overline{\bbF}_q\otimes_k \cH^{(1)}_k= \cH^{(1)}_{\overline{\bbF}_q}$. 
According to \ref{char_k_rational}, a given central character $\eta: Z(G)\rightarrow \overline{\bbF}_q^\times$ is $k$-rational (i.e. takes values in $k^\times$) if and only if $\eta(\pi)\in k^{\times}$. We therefore see that the twisting action restricts to an action of $k$-rational characters on absolutely irreducible supersingular two-dimensional $\cH^{(1)}_{k}$-modules.
\end{Pt}

\begin{Lem}\label{comp_corr_twist}
The correspondence $\ASph(v)^{\rm ss}\rightsquigarrow \rho_{\sL(v)}$ is compatible with twisting by characters.
\end{Lem}

\begin{proof}
Twisting with central characters on semisimple spherical Hecke modules is compatible with the action of $\bbG_m(\bbF_q)\times\bbG_m$ on $S$. The morphism $L$ is $\bbG_m(\bbF_q)\times\bbG_m$-equivariant by theorem \ref{Langlands_morphism}. Under the Galois parametrization \ref{moduli}, the $\bbG_m(\bbF_q)\times\bbG_m$-action on $X_{\bbF_q}$ corresponds to twisting with Galois characters, cf. \ref{qparam_twisting}. Putting all this together, we see that the correspondence $\ASph(v)^{\rm ss}\rightsquigarrow \rho_{\sL(v)}$ is indeed compatible with twisting by characters.
\end{proof}

\begin{Pt} 
Next, we recall the main construction from \cite{GK18} in the case of standard supersingular modules of dimension $n=2$. 

Let $F_{\phi}$ be the special Lubin-Tate group with Frobenius power series $\phi(t)=\pi t + t^q$. Let $F_{\infty}/F$ be the extension generated by all torsion points of $F_{\phi}$ and let $\Gamma = {\rm Gal}(F_{\infty}/F)$. We thus have the category of  \'etale Lubin-Tate $(\varphi,\Gamma)$-modules over $\bbF_q((t))$, cf. \ref{pt_functor}. We identify in the following $\Gamma \simeq o_F^{\times}$ via the Lubin-Tate character $\chi_F$.

To be conform with the notation in \cite[sec. 2.1]{GK18}, we define ${\small
\omega:=u=
\left (\begin{array}{cc}
0 & 1\\
\pi & 0
\end{array} \right)
}$
(this will not lead to confusion with the character of $\bbF_q^{\times}$ denoted by $\omega$ so far).
%in \ref{Fq_action}.) 
In particular $T_{\omega}=U$. Let $M$ be a two-dimensional standard supersingular $\cH^{(1)}_k$-module, arising from a supersingular character $\chi: \cH^{(1)}_{{\rm aff },k}\rightarrow k$ of the affine subalgebra $\cH^{(1)}_{{\rm aff },k}\subset \cH^{(1)}_k$.
Let $e_0\in M$ such that $\cH^{(1)}_{{\rm aff },k}$ acts on $e_0$ via $\chi$ and put $e_1=T^{-1}_{\omega} e_0$. The character $\chi$ determines two numbers 
$0\leq k_0,k_1 \leq q-1$ with $(k_0,k_1)\neq (0,0), (q-1,q-1)$, cf. \cite[Lem. 5.1]{GK18}. One considers $M$ as a $k[[t]]$-module with $t=0$ on $M$. Let $\Gamma = o_F^{\times}$ 
act on $M$ via $$\gamma (m) = T_{e^*(\overline{\gamma})}^{-1} (m)$$ for $\gamma\in o_F^{\times}$ with reduction $\overline{\gamma}\in\bbF_q^{\times}$ and (since $n=2$)
$e^*(\overline{\gamma})=\diag(\overline{\gamma},1) \in\bbT$, cf. \cite[beginning of sec. 4]{GK18}. Moreover, there is a certain
$k[[t]][\varphi]$-submodule $\nabla(M)$ of $$k[[t]][\varphi,\Gamma ]\otimes_{k[[t]][\Gamma]} M \simeq k[[t]][\varphi]\otimes_{k[[t]]} M.$$ The module $\nabla(M)$ is stable under the $\Gamma$-action \cite[Lem. 4.2]{GK18} and thus the quotient
$$\Delta(M) :=  \big( k[[t]][\varphi]\otimes_{k[[t]]} M \big)  / \nabla(M)$$
defines a $k[[t]][\varphi,\Gamma ]$-module. It is torsion standard cyclic with weights $(k_0,k_1)$ in the sense of \cite[sec. 1.3]{GK18}, according to 
\cite[Lemma 5.1]{GK18}.
Let $\Delta(M) ^{*} = \Hom_k (\Delta(M),k)$ be its $k$-linear dual. By a general construction,
the $k((t))$-vector space 
$$\Delta(M) ^{*} \otimes_{k[[t]]} k((t))$$
is then in a natural way an \' etale Lubin-Tate $(\varphi,\Gamma)$-module of dimension $2$. The correspondence $$M\rightsquigarrow  \Delta(M) ^{*} \otimes_{k[[t]]} k((t))$$ extends to a fully faithful functor from 
a certain category of supersingular $\cH^{(1)}_k$-modules to the category of \' etale $(\varphi,\Gamma)$-modules over $k((t))$. 
We write $$V(M):=\sV ( \Delta(M) ^{*} \otimes_{k[[t]]} k((t)) )$$ for 
its composition with the functor $\sV$, cf. \ref{pt_functor}. According to \cite[Cor. 5.5] {GK18}, the map 
$$V\mapsto V(M)$$ 
induces a bijection between (isomorphism classes of) $2$-dimensional supersingular absolutely irreducible $\cH^{(1)}_{k}$-modules and absolutely irreducible representations
${\rm Gal}(\overline{F}/F)\ra {\bf GL_2}(k)$.
\end{Pt}

%Recall that $\varphi^{-1}\in{\rm Gal}(\overline{F}/ F)$ denotes a geometric Frobenius.

\begin{Prop}\label{lem_tensor} One has $V(M\otimes \eta) =  V(M) \otimes \eta $ for any absolutely irreducible supersingular two-dimensional $\cH^{(1)}_{k}$-module $M$ and any character $\eta: {\rm Gal}(\overline{F}/F) \rightarrow k^\times.$ Moreover, if $U^2$ acts by $z_2\in k^\times$ on $M$, then $\det(\varphi^{-1})=z_2$ on $V(M)$. 
\end{Prop}

\begin{proof} The functor $\sV$ respects the tensor product. In our situation, this concretely means the following. 
Write $\eta = \omega_f^s \mu_\lambda$ for a scalar $\lambda\in k^\times$ and $0\leq s \leq q-2.$ Let $D$ be an \' etale $(\varphi,\Gamma)$-module over $k((t))$. We write $D\otimes\eta$ for the $(\varphi,\Gamma)$-module 
equal to the tensor product $D$ by the $1$-dimensional module corresponding to $\eta$: the $\varphi$-action becomes multiplied by the scalar $\lambda$ and the $\Gamma$-action becomes twisted by the character
by $\omega_f^s |_{\Gamma}$, cf. \ref{lem_twist}. Then $\sV(D \otimes\eta) = \sV(D)\otimes \eta$ according to \ref{cor_twist}.
For the first statement, it suffices therefore to check that the functor $M\rightsquigarrow  D(M):=\Delta(M) ^{*} \otimes_{k[[t]]} k((t))$ respects the tensor product with $\eta$.
Let the $\bbT$-action on $M$ be given by $\gamma\in\bbT^\vee/W$.
The corresponding action on $M\otimes\eta$ is then given by 
$\gamma.(\omega_f^s | _{{\bbF}_q^\times})$.
By construction, the $\Gamma$-action on $\Delta(M\otimes\eta)$ is given by the action 
$$a\mapsto T_{e^*(\overline{a})}^{-1} $$ 
on the Hecke module $M\otimes \eta$,
for $a\in o_F^{\times}=\Gamma $ with reduction $\overline{a}\in\bbF_q^{\times}$ and 
$e^*(\overline{a})=\diag(\overline{a},1) \in\bbT$. Since the $\bbT$-action on $M\otimes\eta$ equals $\gamma.(\omega_f^s | _{{\bbF}_q^\times})$ and $$\gamma.(\omega_f^s | _{{\bbF}_q^\times})(e^*(\overline{a}))=\gamma(e^*(\overline{a}))
\omega_f^s(\overline{a})$$
the $\Gamma$-action on $\Delta(M\otimes\eta)$ becomes therefore twisted by the character $\omega_f^{-s} | _{\Gamma}$.
The contragredient action on the dual $\Delta(M\otimes\eta)^*$ and, hence, on $D(M\otimes\eta)$ becomes then twisted by $\omega_f^{s} | _{\Gamma}$, as desired. By construction, the $t^{k_j}\varphi$-action on $\Delta(M\otimes\eta)$ is given by the $T^{-1}_{\omega}$-action on the Hecke module $M\otimes\eta$. Since $U=T_\omega$, this $t^{k_j}\varphi$-action becomes therefore multiplied by $\eta(\pi^{-1})=\lambda^{-1}$. On the dual module $D(M\otimes\eta)$ (where $t$ becomes invertible), the $\varphi$-action becomes therefore multiplied by $\lambda$. We have shown that $D(M\otimes\eta)=D(M)\otimes\eta$, which concludes the proof of the first statement. 

For the second statement, suppose that $U^2$ acts by $z_2\in k^\times$ on $M$. The remarks right before \cite[Thm. 8.8]{GK16} and \cite[Cor. 5.5]{GK18} in the case $d=1$ show that the determinant of geometric Frobenius on $V(M)$ acts by $b$ where $b^{-1}$ equals the $T_{\omega}^{-2}$-action on $M$. Since $U=T_{\omega}$, this implies $b=z_2$ and hence $\det(\varphi^{-1})=z_2$ on $V(M)$. 
\end{proof}

Recall the classification of irreducible $2$-dimensional Galois representations, cf. \ref{irred_n2}. 
\begin{Prop} \label{lem_special}
Let $M$ be a simple supersingular module such that $V(M)\simeq \ind (\omega_{2f}^h)$ where $1\leq h \leq q-1$. Then $M$ has trivial $U^2$-action 
and its $\bbT$-action is given by (the $W$-orbit of) 
the character $\diag(a,b)\mapsto a^{h-1}$.
\end{Prop} 

\begin{proof}
Let $\rho=\ind (\omega_{2f}^h)$. According to \ref{lem_tensor}, the $U^2$-action on $M$ is given by the scalar $\det \rho(\varphi^{-1}) \; =
\omega_f^h(\varphi^{-1})=1.$ To determine the $\bbT$-action, 
let 
$D$ be the \'etale $(\varphi,\Gamma)$-module with $\sV(D)\simeq \rho$, so that
$$\Delta(M) ^{*} \otimes_{k[[t]]} k((t)) \simeq D.$$ 

Now, we make use of the case $n=2$ of the main result \ref{thm_main} of the appendix. It follows that $D$ admits a basis $\{g_0,g_1\}$
such that $$\gamma(g_j)=\overline{f}_{\gamma}(t)^{hq^j/(q+1)}g_j$$ for all $\gamma\in\Gamma$ and $\varphi(g_0)=g_{1}$ and $\varphi(g_1)=-t^{-h(q-1)}g_0$.
Here $$\overline{f}_{\gamma}(t)=\omega_f(\gamma)t/\gamma(t)\in 1+tk[[t]].$$
In particular, $D$ is standard cyclic in the sense of \cite[sec. 1.4]{GK18} with corresponding $\alpha_j:\Gamma\rightarrow k^\times$ given by $\alpha_j=1$ for $j=1,2$ (since 
$\overline{f}_{\gamma}(t)\equiv  1 \mod t$).
Define the triple $$(k_0,k_1,k_2)=(h-1, q-h, h-1)$$ and let $i_j := q-1 -k_{2-j}$, so that $i_0=i_2=q-h$ and $i_1=2q-h-1$. 
Define the triple $(h_0,h_1,h_2)=(0, i_1, i_0 +i_1 q)$. Note that $h_2=h(q-1).$
Put $f_j=t^{h_j}g_j$ for $j=0,1$ and let $D^{\sharp}\subset D$
be the $k[[t]]$-submodule generated by $\{f_0,f_1\}$. Let $(D^{\sharp})^{\ast}$ be the $k$-linear dual. Define $e_i' \in (D^{\sharp})^{\ast}$ via
$e'_i(f_j)=\delta_{ij}$ and $e'_i =0$ on $t D^{\sharp}$. Using the explicit formulae for the $\psi$-operator on $k((t))$ as described in \cite[Lemma 1.1]{GK18} one may follow the argument of \cite[Lemma 6.4]{GK16} and show that $D^{\sharp}$ is a $\psi$-stable lattice in $D$ and that $\{e'_0,e'_1\}$ is a $k$-basis of the $t$-torsion part of  $(D^{\sharp})^{\ast}$ satisfying

$$ t^{k_1} \varphi (e'_0) = e'_1 \hskip5pt \text{and} \hskip5pt  t^{k_0} \varphi (e'_1) = - e'_0. $$

But according to \cite[1.15]{GK18}, there is only one $\psi$-stable lattice in $\Delta(M) ^{*} \otimes_{k[[t]]} k((t))$, 
namely $\Delta(M)^*$. It follows that $\Delta(M)\simeq (D^{\sharp})^{\ast}$ and so the weights of the torsion standard cyclic $k[[t]][\varphi,\Gamma ]$-module
$\Delta(M)$ (in the sense of the definition in \cite[sec. 1.3]{GK18}) are $(k_0,k_1)$. Moreover,
$e'_0,e'_1$ are a $k$-basis of $M$ and $e'_0$ is an eigenvector for the supersingular character
$\chi: \cH^{(1)}_{{\rm aff },k}\rightarrow k$ giving rise to $M$. From $\alpha_j=1$ we deduce from 
the definition of the $\Gamma$-action on $M$, cf. \cite[beginning of sec. 4]{GK18} that 
$T^{-1}_{e^*(\gamma)}=1$ for all $\gamma\in\Gamma$. Hence if $\lambda\in\bbT^\vee$ is the restriction of $\chi$ to $\bbT$, then $$\lambda \circ e^*(\overline{a})=1$$
for any $\overline{a}\in\bbF_q^\times$.
Finally, \cite[Lemma 4.1]{GK18} shows that $k_0\equiv \epsilon_1 \mod (q-1)$ where 
$\epsilon_1$ is such that $\lambda \circ \alpha^\vee (\gamma)^{-1}=\gamma^{\epsilon_1}$
for any $\gamma\in\Gamma$ and the coroot $\alpha^\vee(x)=\diag(x,x^{-1})$, cf.
\cite[discussion before 2.4]{GK18}. This implies that $$\lambda \circ \alpha^{\vee}(\overline{a})^{-1}=\overline{a}^{h-1}$$ for any $\overline{a}\in\bbF_q^\times$.
Since $\diag(\overline{a},\overline{b})=e^*(\overline{a}\cdot\overline{b})\alpha^{\vee}(\overline{b})^{-1}$ we arrive therefore at 
$$ \lambda (\diag(\overline{a},\overline{b})) =\lambda ( e^*(\overline{a}\cdot\overline{b})\alpha^{\vee}(\overline{b})^{-1}) =\overline{b}^{h-1}.$$
\end{proof} 

\begin{Th} \label{comparison_GK}
The correspondence $ \ASph(v)\rightsquigarrow \rho_{\sL(v)}$, when restricted to simple supersingular modules, coincides with the base change from $k$ to $\overline{\bbF}_q$ of the bijection $M\mapsto V(M)$. 
\end{Th}

\begin{proof}
According to \ref{comp_corr_twist} and \ref{lem_tensor}, the correspondences  $\ASph(v)^{\rm ss}\rightsquigarrow \rho_{\sL(v)}$ and $M\mapsto V(M)$ are compatible with twisting. It therefore suffices to compare the two maps on irreducible Galois representations of the form $\rho:=\ind (\omega_{2f}^h)$, for $1\leq h \leq q-1$, cf. \ref{irred_n2}. Let $M$ be such that $V(M)\simeq \rho$. On the other hand, let $\ASph(v_\rho)$ be the supersingular module corresponding to $\rho$ in the bijection \ref{thm_LT}.
Its $U^2$-action is trivial and its $\bbT$-action is given by the highest weight 
$\hw(F(h-1))$, cf. \ref{comp_weights} below.
According to \ref{lem_special}, these actions coincide with the corresponding actions on $M$.
Since both modules are simple supersingular, there is thus an isomorphism $M \simeq \ASph(v_\rho)$.
\end{proof}

\section{Relation to weights}

\begin{Pt} {\bf Weights}.\label{q-weights}
A {\it weight} is an (isomorphism class of an) irreducible  $\overline{\bbF}_q$-representation of 
the finite group ${\rm GL}_2(\bbF_q)$. For any integer $r\geq 0$ consider the $r$-th symmetric power $\Sym^{r}\overline{\bbF}_q^{\oplus 2}$ of the standard ${\rm GL}_2(\bbF_q)$-representation. We denote by $(x,y)$ for a moment the standard basis of $\overline{\bbF}_q^{\oplus 2}$, so that 
$\Sym^{r}\overline{\bbF}_q^{\oplus 2}=\oplus_{i=0,...,r} \overline{\bbF}_q x^{r-i}y^{i}$. The standard action of ${\rm GL}_2(\bbF_q)$ on $\Sym^{r}\overline{\bbF}_q^{\oplus 2}$
is then given by
$$\left (\begin{array}{cc}
a& b \\
c& d
\end{array} \right)(x^{r-i}y^{i})=(ax+cy)^{r-i}(bx+dy)^{i}$$
where $a,b,c,d\in\bbF_q$ are viewed in $\overline{\bbF}_q$ via our fixed embedding $\bbF_q\subset \overline{\bbF}_q$.
Let $$F(r):=\soc_{{\rm GL}_2(\bbF_q)}\Sym^{r}\overline{\bbF}_q^{\oplus 2}$$
be the socle. The representation $F(r)$ is irreducible and contains the highest weight vector $x^r$. The $q(q-1)$ representations $F(r)\otimes\det^s$ for $0\leq r\leq q-1$ and $0\leq s\leq q-2$ exhaust all weights and $F(r)=\Sym^{r}\overline{\bbF}_q^{\oplus 2}$ for $0\leq r\leq p-1$, cf.
 \cite[19.1]{Hu05}. 
\end{Pt}

\begin{Pt} {\bf Compatibility of $\sL$ with weights.}\label{comp_weights}
Recall that we have ordered and labelled the irreducible components of $S$ by the elements of $\mathbf{\whT}(\bbF_q)$ \ref{ordering_comp_S}. Under the isomorphism $\sS^{(1)}_{\bbF_q}$, we have a corresponding ordering of the irreducible components of $\Spec Z(\cH^{(1)}_{\bbF_q})$ by 
$\bbT^{\vee}$. For $\lambda\in\bbT^{\vee}$, we write $\cC^{\lambda}$ for the corresponding irreducible component. On the other hand, to any pair $(r,s(r))$ in the table \ref{types_basic_even_case} or \ref{types_basic_odd_case}, we can associate the weight $F(r)\otimes\det^{s(r)}$. In this way, the irreducible components 
of $X_{\bbF_q}$ can be labelled by ordered pairs of weights $(\sigma,\sigma')$: the irreducible component $\cC^r$ is labelled by the pair $( F(r)\otimes\det^{s(r)}, F(q-3-r)\otimes\det^{s(r)+r+1})$. Finally, we have the highest weight map
\begin{eqnarray*}
\hw:\{\textrm{weights}\} & \lra & \bbT^{\vee}\\
F(r)\otimes\det{}^s & \lmapsto & r(1,0)+s(1,1)  | _{\bbT}.
\end{eqnarray*}
%Finally, to a weight $\sigma$ we also have the (commutative) spherical $\overline{\bbF}_q$-Hecke algebra $\cH^{\sph}(\sigma)$ with coefficients in $\sigma$, e.g. \cite[3.2]{Br07}.
\end{Pt}
\begin{Prop*}\label{comp_weights}
Let $\iota_{\sigma}$ and $\iota_{\sigma'}$ be the embeddings $$
\bbA^1\subset \bbP^1=\cC^{(\sigma,\sigma')}$$ around $0$ and $\infty$ respectively. The morphism $\sL$ induces isomorphisms $$\cC^{\hw(\sigma)}\iso {\rm Im}\; \iota_\sigma\subset\cC^{(\sigma,\sigma')}\hskip10pt \text{and}\hskip10pt \cC^{\hw(\sigma')}\iso {\rm Im} \;\iota_{\sigma'}\subset\cC^{(\sigma,\sigma')}.$$
\end{Prop*}
\begin{proof} The labellings are compatible with the twist by a non-regular character of the form $\omega^s\otimes\omega^s$
and by the determinant character $\det^s$ respectively. It therefore suffices to only consider the basic even case $n=0$ and the basic odd case $n=q-2$. We may also assume $z_2=1$. 

\vskip5pt

Supppose first $n=0$. Then $n=2s$ with $s=0$. The irreducible components of the scheme $\Spec(Z(\cH^{(1)}_{\bbF_q}))_{(\omega^0,1)}$ are labelled by the sequence of pairs of characters $\chi_i, \chi^w_i,$ for $i=0,...,\frac{q-1}{2}$, where $\chi_i:=\omega^i\otimes\omega^{-i}$. By definition of $L$, and hence of $\sL=L\circ \sS^{(1)}_{\bbF_q}$,
$$\sL( \cC^{\chi_i})\subset \cC_{i} \hskip10pt \text{and} \hskip10pt \sL( \cC^{\chi^w_{i+1}})\subset \cC_{i}$$
(the latter if $i<\frac{q-1}{2}$). The weight-label of $\cC_i$ is the pair of weights $(F(r)\otimes\det^{s(r)}, F(q-3-r)\otimes\det^{s(r)+r+1})$ 
where $r=2i$ and $s(r)=-\frac{r}{2}=-i$. For their highest weights we find indeed
$$\hw( F(r)\otimes\det{}^{s(r)} ) =\omega^{r+s(r)}\otimes\omega^{s(r)}=\omega^{i}\otimes\omega^{-i}=\chi_i$$ 
and
$$\hw( F(q-3-r)\otimes\det{}^{s(r)+r+1}) =\omega^{q-2+s(r)}\otimes\omega^{s(r)+r+1}=
\omega^{q-1-(i+1)}\otimes\omega^{i+1}=\chi^w_{i+1}.$$

Now suppose $n=q-2$. Then $n=2s-1$ with $s=\frac{q-1}{2}$. The irreducible components of the scheme $\Spec(Z(\cH^{(1)}_{\bbF_q}))_{(\omega^{q-2},1)}$  are labelled by the sequence of pairs of characters 
$$\chi_i \cdot (\omega^s\otimes\omega^s), \hskip5pt \chi^w_i \cdot (\omega^s\otimes\omega^s),$$ 
for $i=1,...,\frac{q-1}{2}$. 
Note that $\chi_i \cdot (\omega^s\otimes\omega^s)=\omega^{i-1} \otimes\omega^{-i}=:\tilde{\chi}_i$. 
By definition of $L$, 
$$\sL( \cC^{\tilde{\chi}_i})\subset \cC_{i} \hskip5pt (\text{if } i\neq 1) \hskip10pt \text{and} \hskip10pt \sL( \cC^{\tilde{\chi}^w_{i+1}})\subset \cC_{i} \hskip5pt (\text{if } i+1\neq \frac{q-1}{2}).$$
The weight-label of $\cC_i$ is the pair of weights $(F(r)\otimes\det^{s(r)}, F(q-3-r)\otimes\det^{s(r)+r+1})$ 
where $r=2i-1$ and $s(r)=-\frac{r+1}{2}=-i$. For their highest weights we find indeed 
$$\hw( F(r)\otimes\det{}^{s(r)} ) =\omega^{r+s(r)}\otimes\omega^{s(r)}=\omega^{i-1} \otimes\omega^{-i}=\tilde{\chi}_i$$ 
and
$$\hw( F(q-3-r)\otimes\det{}^{s(r)+r+1}) =\omega^{q-2+s(r)}\otimes\omega^{s(r)+r+1}=
\omega^{q-1-(i+1)}\otimes\omega^{(i+1)-1}=\tilde{\chi}^w_{i+1}.$$
It remains to check the cases $i=1$ and $i=\frac{q-1}{2}$, where $L$ is given by the map $t\mapsto t+t^{-1}$. 
In the case $i=1$ the variable $t$ stands for the variable $y$ and by definition of $L$, 
$$ \sL( \cC^{\tilde{\chi}^w_1})\subset \cC_0.$$
The component $\cC_0$ has as weight-label the {\it single} weight $F(q-2)$ and its highest weight is indeed $\hw(F(q-2))=\tilde{\chi}^w_1$. 
In the case $i=\frac{q-1}{2}$ the variable $t$ stands for the variable $x$ and by definition of $L$, 
$$ \sL( \cC^{\chi_{\frac{q-1}{2}}})\subset \cC_{\frac{q-1}{2}}.$$
The component $\cC_{\frac{q-1}{2}}$ has as weight-label the {\it single} weight $F(q-2)\otimes\det^{\frac{q-1}{2}}$ and its highest weight is 
indeed $\hw(F(q-2)\otimes\det^{\frac{q-1}{2}})=\chi_{\frac{q-1}{2}}$. This concludes the proof. 
\end{proof} 

\begin{Pt}
For completeness, we end this section by comparing the above theory, \emph{relative} to $\bbF_q$, with the \emph{absolute} theory over $\bbF_p$.
\end{Pt}

\begin{Pt*}  \label{Breuil} \textbf{$(\bbF_p)^f$-parametrization.}
For $0\leq j\leq f-1$, we denote by $(\Sym^{r}\overline{\bbF}_q^{\oplus 2})^{[j]}$ the $\overline{\bbF}_q$-vector space 
$\oplus_{i=0,...,r} \overline{\bbF}_q x^{r-i}y^{i}$ with the action of ${\rm GL}_2(\bbF_q)$ given by: 

$$
\left (\begin{array}{cc}
a& b \\
c& d
\end{array} \right)(x^{r-i}y^{i})=(a^{p^j}x+c^{p^j}y)^{r-i}(b^{p^j}x+d^{p^j})y^{i},
$$
where $a,b,c,d\in\bbF_q$ are viewed in $\overline{\bbF}_q$ via our fixed embedding $\bbF_q\subset \overline{\bbF}_q$.
Let $r_0,...,r_{f-1}$ and $s$ be integers such that $0\leq r_i\leq p-1$ and $0\leq s\leq q-2$. Then the tensor product representations 

$$ (r_0,...,r_{f-1})\otimes\det{}^s:=(\Sym^{r_0}\overline{\bbF}_q^{\oplus 2})^{[0]}\otimes_{\overline{\bbF}_q} (\Sym^{r_1}\overline{\bbF}_q^{\oplus 2})^{[1]}\otimes_{\overline{\bbF}_q}\cdot\cdot\cdot \otimes_{\overline{\bbF}_q} 
(\Sym^{r_{f-1}}\overline{\bbF}_q^{\oplus 2})^{[f-1]} \otimes\det{}^s$$
form a full set of representatives for the $q^2-q$ weights of ${\rm GL}_2(\bbF_q)$. Indeed, the dictionary with the $q$-parametrization \ref{q-weights} is the following. Let $0\leq r\leq q-1$. Write its $p$-adic expansion $r=r_0 + r_1 p +\cdot\cdot\cdot +r_{f-1}p^{f-1}$ with all $0\leq r_i\leq p-1$.
For any $0\leq s\leq q-2$, one has an isomorphism $$F(r)\otimes\det{}^s \simeq (r_0,...,r_{f-1})\otimes\det{}^s$$
as ${\rm GL}_2(\bbF_q)$-representations, cf. \cite[2.7/2.11/19.1]{Hu05}. 
\end{Pt*}

\begin{Pt*} \textbf{Relation to Serre weights.} In \ref{types_basic_even_case} and  \ref{types_basic_odd_case}
we have associated to a  "horizontal" pair
$$(r,s(r)) \hskip10pt |\hskip10pt    (q-3-r,s(r)+r+1)$$
the niveau $1$ tame inertial type 
$$ \tau:=\left(\begin{array}{cc}
\omega_{f}^{r+1} & 0\\
0 & 1
\end{array} \right)\otimes\omega_{f}^{s(r)}\simeq \left(\begin{array}{cc}
\omega_{f}^{q-2-r} & 0\\
0 & 1
\end{array} \right)
\otimes\omega_{f}^{s(r)+r+1}.
$$
If sufficiently generic, such a type $\tau$ has an associated set $W(\tau)$ of Serre weights\footnote{denoted by $\cD(\tau)$ in \cite{Br07} and called {\it the set of Diamond weights}}, cf. \cite[12.3]{Br07}.
On the other hand, in \ref{comp_weights}, we have associated to such a "horizontal" pair
the two ${\rm GL}_2(\bbF_q)$-weights 
$$(F(r)\otimes{\det}^{s(r)}, F(q-3-r)\otimes{\det}^{s(r)+r+1}).$$ 
These two weights are {\it not} in $W(\tau)$ for generic $\tau$ whenever $f>1$. 

\vskip5pt 

Similarly, in \ref{types_basic_even_case} and  \ref{types_basic_odd_case}
we have associated to a "diagonal" pair
$$(r,s(r)) \hskip10pt \diagup \hskip10pt    (q-1-r,s(r)+r-1)$$
the niveau $2$ tame inertial type 
$$ \tau:=\left(\begin{array}{cc}
\omega_{2f}^{r+1} & 0\\
0 & \omega_{2f}^{q(r+1)}
\end{array} \right)\otimes\omega_{f}^{s(r)}\simeq \left(\begin{array}{cc}
\omega_{2f}^{q-r} & 0\\
0 & \omega_{2f}^{q(q-r)} 
\end{array} \right)
\otimes\omega_{2f}^{s(r-2)+r-1}.
$$
If sufficiently generic, such a type $\tau$ also has an associated set of Serre weights $W(\tau)$, cf.  \cite[8.2]{Br07}. On the other hand, such a "diagonal pair" gives rise to the two ${\rm GL}_2(\bbF_q)$-weights 
$$(F(r)\otimes{\det}^{s(r)}, F(q-1-r)\otimes{\det}^{s(r)+r-1}).$$
These two weights are {\it not} in $W(\tau)$ for generic $\tau$ whenever $f>1$.

\vskip5pt 

In the following, we show the latter statement. The analogous statement in the case of niveau $1$ has a similar proof. So assume $f>1$ and let $x_0,...,x_{f-1}$ be $f$ variables. Let $\cD(x_0,...,x_{f-1})$ be the set of all possible $f$-tuples $\lambda:=(\lambda_0(x_0),...,\lambda_{f-1}(x_{f-1}))$, where 
$\lambda_i(x_i)\in \bbZ\pm x_i$, satisfying the five conditions (i)-(v) listed in \cite[8.1, p.58]{Br07}. 
For example, the first two conditions are:
\begin{itemize}
\item[(i)] $\lambda_0(x_0)\in \{x_0,x_0-1,p-2-x_0,p-1-x_0\}$ and  $\lambda_i(x_i)\in \{x_i,x_i+1,p-2-x_i,p-3-x_i\}$  for $i>0$.
\item[(ii)] if $i>0$ and  $\lambda_i(x_i)\in \{x_i,x_i+1\}$ (or if $i=0$ and $\lambda_0(x_0)\in \{x_0,x_0-1\}$), then $\lambda_{i+1}(x_{i+1})\in \{x_{i+1},p-2-x_{i+1}\}.$ 
\end{itemize} 
One may verify directly that the choice $\lambda_i(x_i)=x_i$ for all $i\geq 0$ meets all five conditions and so 
$$(x_0,...,x_{f-1})\in \cD(x_0,...,x_{f-1}).$$ Similarly, the choice $\lambda(x_0)=p-1-x_0$ and 
$\lambda(x_i)=p-3-x_i$ for $i>0$ meets all five conditions and so $$(p-1-x_0,p-3-x_1,...,p-3-x_{f-1})\in \cD(x_0,...,x_{f-1}).$$
In fact, $\cD(x_0,...,x_{f-1})$ is stable under the transformation $x_0\mapsto p-1-x_0$ and $x_i\mapsto p-3-x_i$ for $i>0$
(this is the "change of variables", as explained in \cite[8.1, p. 59]{Br07}). In general, the set $\cD(x_0,...,x_{f-1})$ has cardinality $| \cD(x_0,...,x_{f-1}) | = 2^f$.

\vskip5pt

To any $\lambda\in \cD(x_0,...,x_{f-1})$, one attaches an expression 
$e(\lambda)(x_0,....,x_{f-1})\in \bbZ \oplus \bigoplus_{i=0,...,f-1} \bbZ x_i$. Its precise definition \cite[8.1, p. 58]{Br07} is not necessary here.
For $\lambda = (x_0,...,x_{f-1})$, one has $e(\lambda)=0$.

\vskip5pt 

As in \ref{types_basic_even_case} and  \ref{types_basic_odd_case}, consider 
$$
\tau = \left(\begin{array}{cc}
\omega_{2f}^{r+1} & 0\\
0 & \omega_{2f}^{q(r+1)}
\end{array} \right)
\otimes\omega_f^{s}
$$
with $-1\leq r\leq q-2$ and $0\leq s\leq q-2$. 
Write the $p$-adic expansion of $r+1$ in the form $$r+1=\sum_{i=0,...,f-1} (r_i+1)p^{i}$$ 
with $0\leq r_i+1 \leq p-1$ for all $i$. We suppose that $\tau$ is {\it generic}, which means precisely that 
$1\leq r_0\leq p-2$ and $0\leq r_i\leq p-3$ for $i>0$, cf. \cite[2.11]{Br07}. In this case, the conditions (i)-(v) in \cite[8.1, p.58]{Br07} ensure that for given $\lambda\in \cD(x_0,...,x_{f-1})$, one has
$0\leq \lambda_i(r_i)\leq p-1$ for all $i$. Then $W(\tau)$ is the set of weights (using the $(\bbF_p)^f$-parametrization) 
$$ (\lambda_0(r_0),...,\lambda_{f-1}(r_{f-1}))\otimes \det{}^{e(\lambda)(r_0,...,r_{f-1})+s}$$
for all $\lambda\in \cD(x_0,...,x_{f-1})$; all these weights are distinct and there are $2^f$ of them, cf. \cite[Prop. 8.3]{Br07}.
For example, taking $\lambda = (x_0,...,x_{f-1})$, one sees that the weight
$(r_0,...,r_{f-1})\otimes\det^s$ is in $W(\tau)$.

\vskip5pt

Now we show that $F(r)\otimes\det^{s}\notin W(\tau)$. Untwisting by $\det^s$, this means that there is {\it no} $\lambda\in  \cD(x_0,...,x_{f-1})$ such that
$F(r)$ can be written in the form 
$$ 
(\lambda_0(r_0),...,\lambda_{f-1}(r_{f-1}))\otimes \det{}^{e(\lambda)(r_0,...,r_{f-1})}.
$$
Write $F(r)$ in the $(\bbF_p)^f$-parametrization:
$$F(r) \simeq (a_0,...,a_{f-1})$$
where $r=\sum_{i=0,...,f-1} a_i p^{i}$ with $0\leq a_i \leq p-1$. 
If there was an equality
$$ (\lambda_0(r_0),...,\lambda_{f-1}(r_{f-1}))\otimes \det{}^{e(\lambda)(r_0,...,r_{f-1})}= (a_0,...,a_{f-1})$$
then uniqueness of the $(\bbF_p)^f$-parametrization and the fact that 
$0\leq \lambda_i(r_i)\leq p-1$ implies that  $\det{}^{e(\lambda)(r_0,...,r_{f-1})}=1$.
Moreover,
$$ r+1 = 1+ \sum_{i=0,...,f-1} a_i p^{i} = \sum_{i=0,...,f-1} (r_i+1)p^{i}.$$
If $a_0=p-1$, then $r_0+1=0$, which is excluded. Hence $a_0\leq p-2$ and unicity of the $p$-adic expansion gives 
$a_0=r_0$ and $a_i=r_i+1$ for $i>0$. Hence $$F(r)\simeq (r_0,r_1+1,...,r_{f-1}+1).$$
If there was $\lambda\in \cD(x_0,...,x_{f-1})$ with $(\lambda_0(r_0),...,\lambda_{f-1}(r_{f-1}))= (r_0,r_1+1,...,r_{f-1}+1)$ then, in particular, 
$\lambda_0(r_0)=r_0$ and $\lambda_1(r_1)=r_1+1$. Comparing with condition (i), the only possibilities for $\lambda_0(x_0)$ and $\lambda_1(x_1)$ are
 $\lambda_0(x_0)=x_0$ and $\lambda_1(x_1)=x_1+1$.  {\it However}, this is excluded by condition (ii). 
This shows that $F(r)\otimes\omega_{f}^{s}\notin W(\tau)$ for any $f>1$.

\vskip5pt 

The same reasoning shows that $F(q-1-r)\otimes\omega_{f}^{s+r}\notin W(\tau)$ for any $f>1$: 
we may start using "other" writing

$$\tau \simeq \left(\begin{array}{cc}
\omega_{2f}^{q-r} & 0\\
0 & \omega_{2f}^{q(q-r)}
\end{array} \right)
\otimes\omega_f^{s+r}.
$$
Write the $p$-adic expansion of $q-r$ as $$q-r=\sum_{i=0,...,f-1} (r_i+1)p^{i},$$ 
with "new" elements $r_0,...,r_{f-1}$. These elements satisfy the same conditions as before, i.e. 
$1\leq r_0\leq p-2$ and $0\leq r_i\leq p-3$ for $i>0$, since the genericity assumption depends only on the isomorphism class of $\tau$. By the same reasoning as above, the weight $(r_0,...,r_{f-1})\otimes \omega_f^{s+r}$ is in $W(\tau)$\footnote{In the "old" variables $r_0,...,r_{f-1}$, the 
corresponding $\lambda\in\cD(x_0,...,x_{f-1})$ is given by $(p-1-x_0,p-3-x_1,...,p-3-x_{f-1})$, cf. \cite[8.1, p. 59]{Br07}.}. Now, the fact that $F(q-1-r)\otimes\omega_{f}^{s+r}\notin W(\tau)$ for any $f>1$ follows word for word as above. Namely, write $F(q-1-r)$ in the $(\bbF_p)^f$-parametrization:
$$F(q-1-r) \simeq (a_0,...,a_{f-1})$$
where $q-1-r=\sum_{i=0,...,f-1} a_i p^{i}$ with $0\leq a_i \leq p-1$. 
If there was an equality
$$ (\lambda_0(r_0),...,\lambda_{f-1}(r_{f-1}))\otimes \det{}^{e(\lambda)(r_0,...,r_{f-1})}=(a_0,...,a_{f-1})$$
then uniqueness of the $(\bbF_p)^f$-parametrization and the fact that 
$0\leq \lambda_i(r_i)\leq p-1$ implies that  $\det{}^{e(\lambda)(r_0,...,r_{f-1})}=1$.
Moreover,
$$ q-r= 1+ \sum_{i=0,...,f-1} a_i p^{i} = \sum_{i=0,...,f-1} (r_i+1)p^{i}.$$
If $a_0=p-1$, then $r_0+1=0$, which is excluded. Hence $a_0\leq p-2$ and unicity of the $p$-adic expansion gives 
$a_0=r_0$ and $a_i=r_i+1$ for $i>0$. Hence $$F(q-1-r)\simeq (r_0,r_1+1,...,r_{f-1}+1).$$
If there was $\lambda\in \cD(x_0,...,x_{f-1})$ with $(\lambda_0(r_0),...,\lambda_{f-1}(r_{f-1}))= (r_0,r_1+1,...,r_{f-1}+1)$, then we conclude as above from condition (i) that $\lambda_0(x_0)=x_0$ and $\lambda_1(x_1)=x_1+1$. {\it However}, this is excluded by condition (ii). 
This shows that $F(q-1-r)\otimes\omega_{f}^{s+r}\notin W(\tau)$ for any $f>1$. \end{Pt*}

 \section{Appendix: Irreducible mod $p$ Lubin-Tate $(\varphi,\Gamma)$-modules}\label{sec_irred}
 
Let $F$ denote a finite extension of $\bbQ_p$, with ring of integers $o_F$ and residue field $\bbF_q$. Let $q=p^f$. Let $\pi\in o_F$ be a uniformizer and let $\overline{F}$ be an algebraic closure of $F$. Let $n\geq 1$ be an integer.

\begin{Pt}
Let $F_{\phi}$ be a Lubin-Tate group for $\pi$, with Frobenius power series $\phi(t)\in o_F[[t]]$. The corresponding ring homomorphism $o_F\rightarrow \End(F_{\phi})$ is denoted by $a\mapsto [a](t)= at +...$. In particular, $[\pi](t)=\phi(t)$.
Let $F_{\infty}/F$ be the extension generated by all torsion points of $F_{\phi}$ and let $$H_F:={\rm Gal}(\overline{F}/F_{\infty})\hskip10pt \text{and} \hskip10pt\Gamma:={\rm Gal}(\overline{F}/F)/ H_F = {\rm Gal}(F_{\infty}/F).$$ 
Let $z=(z_j)_{j\geq 0}$ be a $o_F$-generator of the Tate module of $F_{\phi}$. In particular, for $j\geq 0$
$$z_j=[\pi](z_{j+1})\equiv z_{j+1}^q\mod \pi$$ and $N_{F(z_1)/F}(-z_1)=\pi$. This implies
$$z_{j+1}^q=z_j(1+O(\pi)) \text{~for~} j\geq 1 \text{~~and~~}
z_1^{q-1} = -\pi (1+O(z_1)).$$
The Galois action on the generator $z$ is given by a character 
$\chi_F: {\rm Gal}(\overline{F}/F) \rightarrow o_F^{\times}$, which is surjective and has kernel $H_F$. One has $ \chi_F\equiv \omega_f {\rm ~mod~ } \pi$.
\end{Pt}

\begin{Pt} 
We denote by $\bbC_p$ the completion of an algebraic closure of $\bbQ_p$ and choose an embedding $\overline{F}\subseteq \bbC_p$. Recall that the tilt $\bbC_p^{\flat}$ of the perfectoid field $\bbC_p$ is an algebraically closed and perfect complete non-archimedean field of characteristic $p$. Its valuation ring $o_{\bbC^{\flat}_p}$ is given by the projective limit
$\varprojlim_{x\mapsto x^q} o_{\bbC_p}/\pi o_{\bbC_p}$ and its residue field is $\overline{\bbF}_q$. 
There is a unique multiplicative section $$s: \overline{\bbF}_q\longrightarrow o_{\bbC^{\flat}_p}, a \mapsto (\tau(a){\rm ~mod~} \pi, \tau(a^{q^{-1}}){\rm ~mod~} \pi, \tau(a^{q^{-2}}){\rm ~mod~} \pi,...)$$ 
where $\tau$ denotes the Teichm\"uller map $\overline{\bbF}_{q}\rightarrow o_{\bbC_p}$. There is an inclusion $$\bbF_q((t))\stackrel{\subset}{\longrightarrow}  \bbC_p^{\flat}, ~ t\mapsto  (...,z_j {\rm ~mod}~ \pi,...)$$ and one has $\bbC_p^{\flat}= o_{\bbC^{\flat}_p} [1/t]$. The field $\bbC_p^{\flat}$ is endowed with a continuous action of ${\rm Gal}(\overline{F}/F)$ and 
a Frobenius $\varphi_q$, which raises any element to its $q$-th power.  We let $\bbF_q((t))^{\rm sep}$ denote the separable algebraic closure of $\bbF_q((t))$ inside $\bbC_p^{\flat}$. The field $\bbF_q((t))$ and its separable closure $\bbF_q((t))^{\rm sep}$ inherit the Frobenius action and the commuting ${\rm Gal}(\overline{F}/F)$-action from $\bbC_p^{\flat}$ and there is an
isomorphism $$H_F\stackrel{\simeq}\longrightarrow {\rm Gal}(\bbF_q((t))^{\rm sep}/\bbF_q((t))).$$
\end{Pt}
\begin{Pt}\label{pt_functor}
The theory of Lubin-Tate $(\varphi,\Gamma)$-modules and their relation to Galois representations is developed in \cite{KR09} and \cite{Sch17}.  We only need very basic facts of this theory, and mostly only mod $p$. Note that the power series ring $o_F[[t]]$ has a Frobenius endomorphism and a $\Gamma$-action via
$$\varphi(f)(t)=f( [\pi](t) )\hskip10pt \text{and}\hskip10pt (\gamma f)(t)=f([\chi_F(\gamma)](t))$$ for $f(t)\in o_F[[t]]$. Via reduction mod $\pi$, these actions induce a Frobenius action and a $\Gamma$-action on $\bbF_q[[t]]$ and its quotient field $\bbF_q((t))$.
This allows one to introduce an abelian tensor category of \'etale Lubin-Tate $(\varphi,\Gamma)$-modules over $\bbF_q((t))$. 
 It turns out to be canonically equivalent to the category of continuous finite-dimensional $\bbF_q$-representations of ${\rm Gal}(\overline{F}/F)$, cf. \cite[1.6]{KR09}, \cite[3.2.7]{Sch17}. The functor $\sV$ from  $(\varphi,\Gamma)$-modules to Galois representations is given by $$D \rightsquigarrow \sV(D):=(\bbF_q((t))^{\rm sep}\otimes _{\bbF_q((t))} D)^{\varphi=1}$$ 
where ${\rm Gal}(\overline{F}/F)$ acts diagonally (and via its projection to $\Gamma$ on the second factor).
\end{Pt}

\begin{Pt}\label{kstructure}
Let $k \subset \overline{\bbF}_q$ be a finite extension of $\bbF_q$. One can consider a $k$-representation of ${\rm Gal}(\overline{F}/F)$ as an $\bbF_q$-representation with a $k$-linear structure. Similarly, one may introduce $(\varphi,\Gamma)$-modules over $k((t))=k\otimes_{\bbF_q} \bbF_q((t))$, where $k$ has the trivial Frobenius and $\Gamma$-action. The functor $\sV$ then restricts to an equivalence of categories between  
\'etale $(\varphi,\Gamma)$-modules over $k((t))$ and continuous finite-dimensional $k$-representations of ${\rm Gal}(\overline{F}/F)$.
\end{Pt}

\begin{Pt} 
We fix once and for all an element $y\in \bbF_q((t))^{\rm sep}$ such that 
$$ y^{(q^n-1)/(q-1)}=t.$$
For $g\in {\rm Gal}(\overline{F}/F)$, 
 the power series $$f_g(t)=\frac{\chi_F(g)t}{g(t)}=\frac{\chi_F(g)t}{[\chi_F(g)](t)} \in 1+ t o_F [[t]]$$
depends only on the class of $g$ in $\Gamma$. The same is true for its mod $\pi$ reduction $\overline{f}_g(t)=\omega_{f}(g)t/g(t).$ Note also that the formula  $f^s_g(t)$ defines an element of $o_F[[t]]$ for any $s\in\bbZ_p$.
\end{Pt}

\begin{Lem}\label{lem-aux} One has $g(y)=y \;\omega_{nf}^q(g)\;\overline{f}_g^{-\frac{q-1}{q^n-1}}(t)$ in $\bbF_q((t))^{\rm sep}$ for all $g\in {\rm Gal}(\overline{F}/F_{n})$. 
\end{Lem}

\begin{proof}
This is a generalization of the case $F=\bbQ_p$ treated in \cite[Lem. 2.1.3]{Be10}. Let $j\geq 1$ and choose
$\pi_{nf,j}\in o_{\bbC_p}$ such that $$\pi_{nf,j}^{\frac{q^n-1}{q-1}}=z_j.$$
We write $\pi_j$ for $\pi_{nf,j}$ in the following calculations. Let $g\in {\rm Gal}(\overline{F}/F_{n})$. 
Then $$(g(\pi_j)/\pi_j)^{\frac{q^n-1}{q-1}}=g(z_j)/z_j=\chi_F(g)f^{-1}_g(z_j)$$ and so the quotient of 
$g(\pi_j)/\pi_j$ by $f_g^{-\frac{q-1}{q^n-1}}(z_j)$ is a certain ${\frac{q^n-1}{q-1}}$-th root
of $\chi_F(g)$. Since exponentiation with ${\frac{q^n-1}{q-1}}\in\bbZ_p^{\times}$ is surjective on the subgroup $1+(\pi) \subset o_F^{\times}$
we may write this root as $\tau(\omega_{nf,j}(g))$, with an element $\omega_{nf,j}(g)\in \bbF_{q^n}^\times$, and arrive at 
$$ g(\pi_j)/\pi_j=\tau(\omega_{nf,j}(g))f_g^{-\frac{q-1}{q^n-1}}(z_j).$$

The map $g\mapsto \omega_{nf,j}(g)$ is a character of the group
 ${\rm Gal}(\overline{F}/F_{n})$, since  $$\omega_{nf,j}(g) \equiv g(\pi_j)/\pi_j \mod \mathfrak{m}_{\bbC_p}$$ in the field $\overline{\bbF}_q= o_{\bbC_p}/ \mathfrak{m}_{\bbC_p}$ and this element is fixed by ${\rm Gal}(\overline{F}/F_{n})$.
 Moreover, this character does not depend on the choice of $\pi_j$: a different choice $\pi'_j$ differs from $\pi_j$ by a ${\frac{q^n-1}{q-1}}$-th root of unity, i.e. by an element of $F_{n}$.
 Hence $g(\pi'_j)/\pi'_j=g(\pi_j)/\pi_j$. By this independence, we see (using the element $\pi_{j+1}^q$ as an alternative choice for $\pi_j$) that $$\omega_{nf,j+1}^q=\omega_{nf,j}~ \text{for} ~ j\geq 1.$$
 Moreover, $\pi_{nf,1}^{q^n-1}=z_1^{q-1} = -\pi (1+O(z_1))$ and so $(\pi_{nf}/ \pi_{nf,1})^{q^n-1} \equiv 1 \mod \mathfrak{m}_{\bbC_p}$. The quotient $\pi_{nf}/ \pi_{nf,1} \mod \mathfrak{m}_{\bbC_p}$ is therefore fixed by  ${\rm Gal}(\overline{F}/F_{n})$, in other words $$g(\pi_{nf,1})/\pi_{nf,1} \equiv g(\pi_{nf})/\pi_{nf} \mod \mathfrak{m}_{\bbC_p}$$ for all  $g \in {\rm Gal}(\overline{F}/F_{n})$ and so
 $$  \omega_{nf,1}=\omega_{nf}.$$
 
Now recall that there is an isomorphism $ \varprojlim_{x\mapsto x^q} o_{\bbC_p}\simeq o_{\bbC_p^{\flat}}$ of multiplicative monoids 
given by reduction modulo $\pi$. We use the notation $u=(u^{(j)})$ for elements in the projective limit $\varprojlim_{x\mapsto x^q} o_{\bbC_p}$. The element 
$y\in o_{\bbC_p^{\flat}}$ is given by  $(...,\pi_j {\rm ~mod}~ \pi o_{\bbC_p},...)$. Its preimage
$(y^{(j)})$ under the above isomorphism is therefore given by 
$y^{(j)}=\lim_{m\rightarrow \infty} \pi_{j+m}^{q^m}$. 
By the same argument, the preimage of the element $\overline{f}_g^{-\frac{q-1}{q^n-1}}(t)$ has coordinates

$$\overline{f}_g^{-\frac{q-1}{q^n-1}}(t)^{(j)}=\lim_{m\rightarrow \infty} (f_g^{-\frac{q-1}{q^n-1}}(z_{j+m}))^{q^m}.$$

The composite map $s: \overline{\bbF}_{q}\rightarrow  o_{\bbC_p^{\flat}} \simeq \varprojlim_{x\mapsto x^q} o_{\bbC_p}$, which we also denote by $s$, is given by
$a \mapsto (\tau(a), \tau(a^{q^{-1}}), \tau(a^{q^{-2}}),...)$.
Since $$s(\omega_{nf}(g)^q)^{(j)}=\tau(\omega_{nf}(g)^{q^{-j+1}})=\tau(\omega_{nf,j}(g)),$$ we may put everything together and obtain
$$ \frac{g(y^{(j)})}{y^{(j)}}=\lim_{m\rightarrow \infty} (\frac{g(\pi_{j+m})}{\pi_{j+m}})^{q^m}=\tau( \omega_{nf,j}(g))\lim_{m\rightarrow \infty} (f_g^{-\frac{q-1}{q^n-1}}(z_{j+m}))^{q^m}=s(\omega_{nf}(g)^q)^{(j)}
\overline{f}_g^{-\frac{q-1}{q^n-1}}(t)^{(j)}.$$ Reducing this equation modulo $\pi$ yields the assertion of the lemma. 
\end{proof}

We now consider the
 $(\varphi,\Gamma)$-modules associated to the
irreducible Galois representations of the form $\ind(\omega_{nf}^h)$. 

\begin{Th}\label{thm_main} The \'etale Lubin-Tate $(\varphi,\Gamma)$-module associated to an irreducible Galois representation of the form 
$\ind(\omega_{nf}^h)$
is defined over the ring $\bbF_q((t))$ and admits a basis $e_0,e_1,...,e_{n-1}$ in which 
$$\gamma(e_j)=\overline{f}_{\gamma}(t)^{hq^j(q-1)/(q^n-1)}e_j$$ for all $\gamma\in\Gamma$. Moreover, one has $\varphi(e_j)=e_{j+1}$ and $\varphi(e_{n-1})=(-1)^{n-1}t^{-h(q-1)}e_0$.
\end{Th}

\begin{proof}
Let $D$ be the $(\varphi,\Gamma)$-module described in the statement and let $W=\sV(D)$. With $x=t^h e_0 \wedge ... \wedge e_{n-1}$, one has
$$ \varphi(x)= \varphi(t)^h (-1)^{n-1}t^{-h(q-1)}  e_1 \wedge ... \wedge e_{n-1} \wedge e_0 = t^{qh-h(q-1)} e_0 \wedge ... \wedge e_{n-1} = x.$$
Moreover, 
$$ \gamma(t)^h \prod_{ j=0}^{n-1} \overline{f}_{\gamma}^{h q^j (q-1)/(q^n-1)}(t) = 
( \omega_f(\gamma) t / \overline{f}_{\gamma}(t) ) ^h \overline{f}_{\gamma}^{h(q-1)/(q^n-1) \sum_{ j=0}^{n-1} q^j}  = \omega_f(\gamma)^h t^h$$
which implies $\gamma(x)= \omega_f(\gamma)^h x$ for all $\gamma\in\Gamma$. So $\det W = \omega_f ^h$. Put $k:=\bbF_{q^n}$ as a coefficient field, i.e. endowed with the trivial Frobenius action. To complete the proof, it remains to check that the restriction of $k\otimes_{\bbF_q} W$ to the inertia subgroup $\cI$ is given by 
$\omega_{nf}^{h} \oplus \omega_{nf}^{qh}\oplus ...\oplus \omega_{nf}^{q^{n-1}h}$. There is a ring isomorphism  
$$k \otimes_{\bbF_q}\bbF_q((t))^{\rm sep}\stackrel {\simeq}{\longrightarrow}\prod_{j=0}^{n-1}\bbF_q((t))^{\rm sep},\; x\otimes z \mapsto (\varphi_q^j (x)z)$$ where $\varphi_q$ is the $q$-Frobenius on $k$.
The induced Frobenius and ${\rm Gal}(\overline{F}/F_{n})$-action on $\prod_{j=0}^{n-1} \bbF_q((t))^{\rm sep} $ are given as 
$$
\begin{array}{ccl}
 \varphi ((x_0,...,x_{n-1})) & =&  (\varphi_q (x_{n-1} ), \varphi_q(x_0),...,\varphi_q (x_{n-2})) \\
&&\\
g ((x_0,...,x_{n-1}))&=  &(g (x_0),...,g( x_{n-1} )). 
\end{array}
$$
Choose $\alpha\in\overline{\bbF}_q \subset \bbF_q((t))^{\rm sep} $ such that $\alpha^{q^n-1}=(-1)^{n-1}$ and define the elements 
$$
\begin{array}{ccl}
 v_0 & =&  (\alpha y^h,0,...,0)e_0+ (0, \alpha^q y^{qh},0,...,0)e_1+...+  (0,...,0, \alpha^{q^{n-1}} y^{q^{n-1}h})e_{n-1}\\
&&\\
v_1 &=  & (0,\alpha y^h,0,...,0)e_0+ (0,0, \alpha^q y^{qh},0,...,0)e_1+...+  (\alpha^{q^{n-1}} y^{q^{n-1}h},0,...,0)e_{n-1} \\
\vdots &&\\
v_{n-1} &=  & (0,...0,\alpha y^h)e_0+ (\alpha^q y^{qh},0,...,0)e_1+...+  (0,...,\alpha^{q^{n-1}} y^{q^{n-1}h},0)e_{n-1} . 
\end{array}
$$
By definition of $D$, the vectors $e_i$ form a $\bbF_q((t))$-basis for $D$ and it follows easily 
that the vectors $v_i$ form a $k \otimes_{\bbF_q}  \bbF_q((t))^{\rm sep} $-basis
 for $k\otimes_{\bbF_q} ( \bbF_q((t))^{\rm sep} \otimes_{\bbF_q((t))}  D )$. Moreover, 
 
 $$ 
 \begin{array}{ccl}
 
 \varphi (  (0,...,0, \alpha^{q^{n-1}} y^{q^{n-1}h})e_{n-1}) &= & ( \alpha^{q^{n}} y^{q^{n}h},0,...,0)  \varphi ( e_{n-1}) \\
 &&\\
&=&  ( \alpha^{q^{n}} y^{q^{n}h},0,...,0)) (-1)^{n-1}t^{-h(q-1)} e_{0} =  (\alpha y^h,0,...,0)e_0
 \end{array}
 $$
 since  $\alpha^{q^{n}} = (-1)^{n-1} \alpha$ and $y^{q^n} t^{1-q}=y$. This means 
  $$
\begin{array}{ccl}
 \varphi (v_0) & =&  (0,\alpha^q y^{qh},0,...,0) \varphi (e_0)+ ( 0,0, \alpha^{q^2} y^{q^2h},0,...,0) \varphi (e_1)+...+ ( \alpha^{q^{n}} y^{q^{n}h},0,...,0) \varphi (e_{n-1}) \\
&&\\
&=  &   (0,\alpha^q y^{qh},0,...,0) e_1+ ( 0,0, \alpha^{q^2} y^{q^2h},0,...,0) e_2+...+  (\alpha y^h,0,...,0)e_0 \\
&=& v_0.
\end{array}
$$
Similarly, one shows $\varphi(v_j)=v_j$ for $j\geq 1$, so that $$v_0,...,v_{n-1} \in k \otimes_{\bbF_q} ( \bbF_q((t))^{\rm sep}  \otimes_{\bbF_q((t))}  D )^{\varphi=1}= k \otimes_{\bbF_q} \sV(D)= k \otimes_{\bbF_q}  W.$$
 Now if $g\in{\rm Gal}(\overline{F}/F_{n})$, then $g(y)=y \omega_{nf}^q(g)c_g$ with $c_g:=\overline{f}_g^{-\frac{q-1}{q^n-1}}(t)$ by lemma \ref{lem-aux} and
 $ g(e_j)=c_g ^{-q^j h}e_j$ by definition of $D$. Hence
 $$ g(y)^{q^jh}g(e_{j})= (y \omega_{nf}^q(g) )^{q^jh}e_j.$$ 
 If $g\in\cI$, then $g(\alpha)=\alpha$ and then altogether
 
 $$
\begin{array}{ccl}
 g (v_0) & =&  (\alpha g(y)^{h},0,...) g (e_0)+ ( 0,\alpha^{q} g(y)^{qh},0,...) g (e_1)+...+ (0,...,\alpha^{q^{n-1}} g(y)^{q^{n-1}h}) g (e_{n-1}) \\
&&\\
&=  & \omega_{nf}^{qh}(g)\cdot (  (\alpha y^{h},0,...) e_0 + ( 0,\alpha^{q} y^{qh},0,...) e_1+...+ (0,...,\alpha^{q^{n-1}} y ^{q^{n-1}h}) e_{n-1})   \\
&&\\
&=&  \omega_{nf}^{qh}(g)\cdot v_0,
\end{array}
$$
where $\cdot$ refers to the left $k$-structure of $\prod_{j=0}^{n-1} \bbF_q((t))^{\rm sep}$.
Similarly, one shows $g(v_j)=  \omega_{nf}^{q^{1-j}h}(g)v_j$ for all $j\geq 1$ and $g\in\cI$. Since $\omega_{nf}^{q^n}= \omega_{nf}$ and hence $\omega_{nf}^{q^{1-j}h}= \omega_{nf}^{q^{n+1-j}h}$, this proves that 
the restriction of $k \otimes_{\bbF_q} W$ to $\cI$ is given by the sum of the characters
$\omega_{nf}^{h}, \omega_{nf}^{qh},...,\omega_{nf}^{q^{n-1}h}$. 
\end{proof}
\begin{Pt}\label{char_k_rational}
One may pass from irreducible representations of the form $\ind(\omega_{nf}^h)$ to general irreducible representations by twisting with characters, cf. \ref{irred_k_reps}.
Note that any character $\eta: {\rm Gal}(\overline{F}/F)\rightarrow \overline{\bbF}^{\times}_q$ can be written in the form 
 $\omega_f^s \mu_{\lambda}$ for a scalar $\lambda\in\overline{\bbF}^{\times}_q $ and $0\leq s \leq q-2.$ In particular, $\eta$ is $k$-rational for a finite extension $k\subset \overline{\bbF}_q$ of $\bbF_q$ (i.e. $\eta$ takes values in $k$) if and only if $\eta(\varphi)\in k^{\times}$.
 \end{Pt}
 
\begin{Lem}\label{lem_twist} 
Let $k\subset \overline{\bbF}_q$ be a finite extension of $\bbF_q$. The $(\varphi,\Gamma)$-module associated to a Galois character of the form 
$\omega_f^s \mu_{\lambda}$ with $\lambda\in k^{\times}$ admits a basis $e$ such that $\varphi(e) = \lambda \cdot e$ and $\gamma (e) = \omega_f^s (\gamma)\cdot e$ for all $\gamma\in\Gamma$.
\end{Lem}

\begin{proof}
Since the functor $\sV$ preserves the tensor product, we may discuss the two characters  $\omega_f^s$ and $\mu_{\lambda}$ separately. 
For the twists by a character of $\Gamma$, such as $\omega_f^s$, see \cite[Remark 4.6]{SV16}.
So let $V=\mu_{\lambda}=k$ and let 
$$D(V)=(\bbF_q((t))^{\rm sep}\otimes _{\bbF_q} V)^{H_F}  $$ be the associated $(\varphi,\Gamma)$-module. 
It is instructive to check the case $k=\bbF_{q}$ first. Here, we choose $\beta\in \overline{\bbF}_{q}$ with $\beta^{q-1}=\lambda$ and 
put $e= \beta \otimes 1\in  \bbF_q((t))^{\rm sep} \otimes_{\bbF_q}  V$. Since $\beta \neq 0$, we have $e\neq 0$. Moreover, $\cI$ acts trivial on $e$ and for $\varphi\in {\rm Gal}(\overline{F}/F)$
we find 
$$ \varphi ( e ) = \varphi ( \beta ) \otimes \varphi ( 1 ) = \beta^q \otimes \lambda^{-1} = \beta \lambda \otimes \lambda^{-1} = \beta \otimes 1 = e.$$
Hence $e$ is indeed ${\rm Gal}(\overline{F}/F)$-invariant. Moreover, if $\phi$ denotes the Frobenius endomorphism on $D(V)$ we have 
$$\phi (e ) = \phi ( \beta ) \otimes 1 = \beta^q \otimes 1 =  \lambda \beta \otimes 1 = \lambda e.$$
Now suppose that $k=\bbF_{q^n}$ for some $n$ and $\lambda\in k^{\times}$. We use the ring isomorphism
$$k \otimes_{\bbF_q}\bbF_q((t))^{\rm sep}\stackrel {\simeq}{\longrightarrow}\prod_{j=0}^{n-1}\bbF_q((t))^{\rm sep},\; x\otimes z \mapsto (\varphi_q^j (x)z)$$ where $\varphi_q$ is the $q$-Frobenius on $k$. It is ${\rm Gal}(\overline{F}/F_n)$-equivariant, where the Galois action on the right-hand side is componentwise (see proof of the above theorem). By the normal basis theorem, there is $x\in k^{\times}$ such that its conjugates $\varphi_q^j (x)$ are linearly independent over ${\bbF_q}$.
The $j$-th copy $\bbF_q((t))^{\rm sep}$ in the above product has therefore a $\bbF_q((t))^{\rm sep}$-basis element $e_j := \varphi_q^j (x)\in k=V$ on which $\cI$ acts trivial and on which the element $\varphi^n\in {\rm Gal}(\overline{F}/F_n)$ acts by $\lambda^{-n}$. Choose $\beta\in \overline{\bbF}_{q}$ such that $\beta^{q^n-1}=\lambda^{n}$ and put 
$v_j = \beta e_j $. Then $\cI$ obviously acts trivial on $v_j$ and the same holds for $\varphi^n$, since 
$$\varphi^n(v_j) = \varphi^n (\beta) \varphi^n (e_j ) =  \beta^{q^n} \lambda^{-n} e_j  = \beta \lambda^{n}  \lambda^{-n} e_j =v_j.$$
Hence, $\cI$ and $\varphi^n$ act trivial on $(v_j)\in \prod_{j=0}^{n-1}\bbF_q((t))^{\rm sep}$ and then also on its preimage $v=x \otimes \beta\in  k \otimes_{\bbF_q}\bbF_q((t))^{\rm sep}$. 
Note that $v\neq 0$ since $x,\beta \neq 0$. Write $N= \prod_{j=0}^{n-1} \varphi^j$ and $e= N(v)$. Then $e$ is fixed by $\cI$ (since $\cI$ is normalized by the $\varphi^j$) and is fixed by $\varphi$ by construction. Hence, $e$ is ${\rm Gal}(\overline{F}/F)$-invariant. Note that $e\neq 0$, since $e=N(x) \otimes N(\beta)$ and $N(x),N(\beta) \neq 0$ and so $e$ is indeed a basis element of $D(V)$ on which $\Gamma$ acts trivial. Finally, write $e=\sum_{j=0}^{n-1}  \varphi^j_q(x)\otimes z_j $
with $z_j \in \bbF_q((t))^{\rm sep}$. The Frobenius endomorphism $\phi$ on $D(V)$ satisfies
$$\phi ( e ) = \sum_j \varphi^j_q(x)\otimes \varphi (z_j)  = \varphi ( \sum_j  \varphi^{-1} (\varphi^j_q(x)) \otimes z_j )= 
 \varphi ( \sum_j  \lambda\varphi^j_q(x)\otimes z_j ) = \lambda \varphi (e) =\lambda e.$$ 
 \end{proof}

\begin{Cor} \label{cor_twist}
Let $k\subset\overline{\bbF}_q$ be a finite extension of $\bbF_q$. The $(\varphi,\Gamma)$-module associated to an irreducible Galois representation of the form 
$( \ind(\omega_{nf}^h)) \otimes \omega_f^s \mu_{\lambda}$, with $\lambda^n \in k^{\times}$,
is defined over the ring $k((t))$ and admits a basis $e_0,e_1,...,e_{n-1}$ in which 
$$\gamma(e_j)=\omega_f(\gamma)^s \overline{f}_{\gamma}(t)^{hq^j(q-1)/(q^n-1)}e_j$$ for all $\gamma\in\Gamma$. Moreover, one has
$\varphi(e_j)=\lambda e_{j+1}$ and $\varphi(e_{n-1})=(-1)^{n-1}t^{-h(q-1)}\lambda e_0$.
\end{Cor}
\begin{proof}
This follows from the preceding lemma and the theorem. The fact that the module is defined over $k((t))$ comes from \ref{irred_k_reps} and
 \ref{kstructure}.
\end{proof}

\vskip10pt 

\noindent {\small Cédric Pépin, LAGA, Université Paris 13, 99 avenue Jean-Baptiste Clément, 93 430 Villetaneuse, France \newline {\it E-mail address: \url{cpepin@math.univ-paris13.fr}} }

\vskip10pt

\noindent {\small Tobias Schmidt, Univ Rennes, CNRS, IRMAR - UMR 6625, F-35000 Rennes, France \newline {\it E-mail address: \url{tobias.schmidt@univ-rennes1.fr}} }


\begin{thebibliography}{99}

\bibitem[Be10]{Be10} {\sc L. Berger},
{\it On some modular representations of the Borel subgroup 
of ${\rm GL_2}(\bbQ_p)$}, Comp. Math. \textbf{146}(1) (2010), 58-80.

\bibitem[Bo11]{Bo11} {\sc C. Bonnaf\'{e}}, {\it Representations of {${\rm SL}_2(\Bbb F_q)$}}, Algebra and Applications \textbf{13} (Springer-Verlag London, Ltd., London, 2011).





\bibitem[Br07]{Br07} {\sc C. Breuil}, {\it Representations of Galois and of  ${\rm GL_2}$ in characteristic $p$}, Course at Columbia University (2007), \url{https://www.imo.universite-paris-saclay.fr/~breuil/PUBLICATIONS/New-York.pdf}.

\bibitem[Br15]{Br15} {\sc C. Breuil}, {\it Induction parabolique et $(\varphi,\Gamma)$-modules.}, Algebra and Number Theory. \textbf{9} (2015), 2241-2291.


\bibitem[Be11]{Be11} {\sc L. Berger},
{\it La correspondance de Langlands locale $p$-adique pour ${\rm GL_2}(\bbQ_p)$}, Astérisque \textbf{339} (2011), 157-180.

%\bibitem[Br03]{Br03} {\sc C. Breuil}, {\it Sur quelques représentations modulaires et $p$-adiques de $GL_2(\bbQ_p)$. I}, Compositio Math. \textbf{138} (2003), 165-188.

%\bibitem[Br97]{Br97} {\sc M. Brion}, {\it Equivariant Chow groups for torus actions}, Transformations Groups, Vol. 2, Nr. 3, 1997, pages 225-267.

%\bibitem[CG97]{CG97} {\sc N. Chriss and V. Ginzburg},
%{\it Representation theory and complex geometry}, Birkhäuser, Boston,1997.

%\bibitem[C10]{C10} {\sc P. Colmez}, {\it Représentations de $GL_2(\bbQ_p)$ et $(\varphi,\Gamma)$-modules}, Astérisque \textbf{330} (2010), 281-509.

%\bibitem[D73]{D73} {\sc M. Demazure}, {\it Invariants symétriques entiers des groupes de Weyl et torsion}, Invent. Math. 21, pages 287-301, 1973.

%\bibitem[D74]{D74} {\sc M. Demazure}, {\it Désingularisation des variétés de Schubert généralisées}, Ann. Sci. \'Ecole Norm. Sup. (4), tome 7 (1974), 53-88. Collection of articles dedicated to Henri Cartan on the occasion of his 70th birthday, I.

\bibitem[Em19]{Em19} {\sc M. Emerton},  {\it Localizing ${\rm GL_2}(\bbQ_p)$-representations}, Talk at Padova School on Serre conjectures and the $p$-adic Langlands program 2019, \url{https://mediaspace.unipd.it/channel/School+}\newline \url{on+Serre+conjectures+and+the+p-adic+Langlands+program/119214951}.

%\bibitem[EG96]{EG96} {\sc D. Edidin and W. Graham}, {\it Equivariant Intersection Theory}, Invent. Math. 131, pages 595-634, 1996.

\bibitem[GK16]{GK16} {\sc E. Grosse-Klönne}, {\it From pro-$p$-Iwahori-Hecke modules to $(\varphi,\Gamma)$-modules, I}, Duke Math. Journal 165 No. 8 (2016), 1529-1595.

\bibitem[GK18]{GK18} {\sc E. Grosse-Klönne}, {\it Supersingular Hecke modules as Galois representations},  Algebra Number Theory 14 (2020), no. 1, 67–118. 

\bibitem[H09]{H09} {\sc F. Herzig}, {\it The weight in a Serre-type conjecture for tame
$n$-dimensional Galois representations}, Duke Math. J. 149(1), pages 37-116.


\bibitem[Hu05]{Hu05} {\sc J.E. Humphreys}, {\it Modular Representations of Finite Groups of Lie Type}, London Math. Soc. Lecture Notes Series 326, Cambridge University Press, 2005.

%\bibitem[KK86]{KK86} {\sc B. Kostant, S. Kumar}, {\it The Nil Hecke ring and cohomology of $G/P$ for a Kac-Moody group $G$},
%Advances. Math., \textbf{62} (1), 187-237, 1986.

\bibitem[KL87]{KL87} {\sc D. Kazhdan, G. Lusztig}, {\it Proof of the Deligne-Langlands conjecture for Hecke algebras},
Invent. Math., \textbf{87} (1), 153-215, 1987.

%\bibitem[KL87]{KL87} {\sc D. Kazhdan, G. Lusztig}, {\it Proof of the Deligne-Langlands conjecture for Hecke algebras},
%Invent. Math., \textbf{87} (1), 153-215, 1987.


\bibitem[KR09]{KR09} {\sc M. Kisin, W. Ren}, {\it Galois representations and Lubin-Tate groups},
Doc. Math., \textbf{14}, 441-461, 2009.

\bibitem[O09]{O09} {\sc R. Ollivier}, {\it Le foncteur des invariants sous l'action du pro-p-Iwahori de $GL(2,F)$}, 
J. f\"ur die reine und angewandte Mathematik \textbf{635} (2009), 149-185.

%\bibitem[O14]{O14} {\sc R. Ollivier}, {\it Compatibility between Satake and Bernstein isomorphisms in characteristic $p$}, Algebra and Number Theory \textbf{8}(5) (2014), 1071-1111.

\bibitem[PS]{PS} {\sc C. Pépin, T. Schmidt}, {\it Generic and Mod $p$ Kazhdan-Lusztig Theory for $GL_2$}, Preprint (2020) arXiv:2007.01364v1.

\bibitem[PS2]{PS2} {\sc C. Pépin, T. Schmidt}, {\it A semisimple mod $p$ Langlands correspondence in families for $GL_2(\bbQ_p)$}, Preprint 2023, submitted for publication. 
arXiv:2009.07328v2.

\bibitem[Sch17]{Sch17} {\sc P. Schneider}, {\it Galois representations and {$(\varphi,\Gamma)$}-modules}, Cambridge Studies in Advanced Mathematics \textbf{164} (Cambridge University Press, 2017).

\bibitem[SV16]{SV16} {\sc P. Schneider, O. Venjakob}, {\it Coates-{W}iles homomorphisms and {I}wasawa cohomology for
              {L}ubin-{T}ate extensions}, Elliptic curves, modular forms and {I}wasawa theory, Springer Proc. Math. Stat. \textbf{188}, 401-468, 2016.


\bibitem[Se72]{Se72} {\sc J.-P. Serre}, {\it Propri\'{e}t\'{e}s galoisiennes des points d'ordre fini des courbes
   elliptiques}, Invent. Math., \textbf{15}(4) (1972), 259-331.
   




%\bibitem[St75]{St75} {\sc R. Steinberg}, {\it On a theorem of Pittie}, Topology \textbf{14} (1975), 173-117.

%\bibitem[Th87]{Th87} {\sc R.W. Thomason}, {\it Algebraic $K$-theory of group scheme actions}, Ann. of Math. Stud. \textbf{113} (1987), 539-563.

\bibitem[V94]{V94} {\sc M.-F. Vigneras}, {\it A propos d'une conjecture de Langlands modulaire}, in: {\it Finite reductive groups}, Ed. M. Cabanes, Prog. Math., vol. 141 (Birkhäuser, Basel, 1997).

%\bibitem[V96]{V96} {\sc M.-F. Vigneras}, {\it Repr\'esentations $\ell$-modulaires d'un groupe r\'eductif p-adique}, Prog. Math., vol. 137 (Birkhäuser, Basel, 1996).

\bibitem[V04]{V04} {\sc M.-F. Vigneras} {\it Representations modulo $p$ of the $p$-adic group $GL(2,F)$}, Compositio Math. 140 (2004) 333-358.

%\bibitem[V05]{V05} {\sc M.-F. Vigneras} {\it Pro-$p$-Iwahori Hecke ring and supersingular $\overline{\bbF}_p$-representations}, Math. Ann. \textbf{331} (2005), 523-556. + Erratum


%\bibitem[V06]{V06} {\sc M.-F. Vigneras} {\it Algèbres de Hecke affines génériques}, Representation Theory \textbf{10} (2006), 1-20.

%\bibitem[V14]{V14} {\sc M.-F. Vigneras}, {\it The pro-$p$-Iwahori Hecke algebra of a reductive $p$-adic group II}, Compositio Math. Muenster J. Math. \textbf{7} (2014), 363-379. + Erratum

  %  \bibitem[V15]{V15} {\sc M.-F. Vigneras}, {\it The pro-$p$-Iwahori Hecke algebra of a reductive $p$-adic group V (Parabolic induction)}, Pacific J. of Math. \textbf{279} (2015), Issue 1-2, 499-529.

%\bibitem[V16]{V16} {\sc M.-F. Vigneras}, {\it The pro-$p$-Iwahori Hecke algebra of a reductive $p$-adic group I}, Compositio Math. \textbf{152} (2016), 693-753.

%\bibitem[V17]{V17} {\sc M.-F. Vigneras}, {\it The pro-$p$-Iwahori Hecke algebra of a reductive $p$-adic group III (Spherical Hecke algebras and supersingular modules)}, Journal of the Institue of Mathematics of Jussieu \textbf{16} (2017), Issue 3, 571-608. + Erratum

%\bibitem[Yo08]{Yo08} {\sc T. Yoshida}, {\it Local class field theory via {L}ubin-{T}ate theory}, Ann. Fac. Sci. Toulouse Math., \textbf{17}(6) (2008), 411-438.




\end{thebibliography}
\end{document}